\documentclass[11pt]{article}
\usepackage{graphicx,dcolumn,amsmath,latexsym,amssymb}
\usepackage{amsfonts,subfigure,caption}
\usepackage{color,bm}
\usepackage{placeins}
\usepackage{hyperref}

\def\H{{\mathcal{H}}}
\def\O{{\mathcal{O}}}

\def\Q{{\mathcal{Q}}}
\def\R{{\mathcal{R}}}
\def\E{{\mathcal{E}}}
\def\F{{\mathcal{F}}}
\def\K{{\mathcal{K}}}
\def\L{{\mathcal{L}}}

\def\G{{\mathcal{G}}}
\def\T{{\mathcal{T}}}

\def\PP{{\mathbb{P}}}
\def\RR{{\mathbb{R}}}
\def\CC{{\mathbb{C}}}

\def\ZZ{{\mathbb{Z}}}

\def\e{{\mathbf e}}

\def\k{{\mathbf k}}
\def\x{{\bm x}}

\def\u{{\bm u}}

\def\x{{\mathbf x}}

\def\z{{\mathbf z}}
\def\f{{\mathbf f}}

\def\u{{\mathbf u}}
\def\v{{\mathbf v}}

\def\0{{\mathbf 0}}

\def\bnabla{\boldsymbol{\nabla}}
\def\bDelta{\boldsymbol{\Delta}}

\def\Bmp#1{ \begin{minipage}{#1} }
\def\Emp{ \end{minipage} }
\def\Bmpc#1{ \begin{minipage}[c]{#1} }
\def\Bmpt#1{ \begin{minipage}[t]{#1} }
\def\Bmpb#1{ \begin{minipage}[b]{#1} }

\newcommand{\uvec}{\mathbf{u}}

\newcommand{\laplacian}{\Delta}
\newcommand{\rot}{\bnabla\times}
\newcommand{\tu}{\widetilde{\mathbf{u}}}

\newcommand{\tuE}{\widetilde{\mathbf{u}}_{\E_0}}

\newcommand{\tuET}{\widetilde{\uvec}_{0;\E_0,T}}
\newcommand{\tuST}{\widetilde{\uvec}_{0;S,T}}
\newcommand{\tuBT}{\widetilde{\uvec}_{0;B,T}}
\newcommand{\tuKT}{\widetilde{\uvec}_{0;\K_0,T}}

\newcommand{\uST}{\uvec_{0;S,T}}
\newcommand{\uBT}{\uvec_{0;B,T}}

\newcommand{\argmax}{\operatorname{argmax}}

\newcommand{\Id}{\operatorname{Id}}

\newtheorem{problem}{Problem}

\begin{document}

\title{\vspace*{-2.5cm} Searching for Singularities in
    Navier-Stokes Flows Based on the Ladyzhenskaya-Prodi-Serrin Conditions}
\author{Di Kang$^1$ and Bartosz Protas$^{1,}$\thanks{Corresponding Author, Email: {\tt bprotas@mcmaster.ca}}} \date{
  $^1$ Department of Mathematics and Statistics, McMaster University \\
  Hamilton, ON, Canada \\ \medskip \today}
\maketitle

\begin{abstract}
  In this investigation we perform a systematic computational search
  for potential singularities in 3D Navier-Stokes flows based on the
  Ladyzhenskaya-Prodi-Serrin conditions. They assert that if the
  quantity $\int_0^T \| \mathbf{u}(t) \|_{L^q(\Omega)}^p \, dt$, where
  $2/p+3/q \le 1$, $q > 3$, is bounded, then the solution
  $\mathbf{u}(t)$ of the Navier-Stokes system is smooth on the
  interval $[0,T]$. In other words, if a singularity should occur at
  some time $t \in [0,T]$, then this quantity must be unbounded. We
  have probed this condition by studying a family of variational PDE
  optimization problems where initial conditions $\mathbf{u}_0$ are
  sought to maximize $\int_0^T \| \mathbf{u}(t) \|_{L^4(\Omega)}^8 \,
  dt$ for different $T$ subject to suitable constraints. These
  problems are solved numerically using a large-scale adjoint-based
  gradient approach. Even in the flows corresponding to the optimal
  initial conditions determined in this way no evidence has been found
  for singularity formation, which would be manifested by unbounded
  growth of $\| \mathbf{u}(t) \|_{L^4(\Omega)}$.  However, the maximum
  enstrophy attained in these extreme flows scales in proportion to
  $\mathcal{E}_0^{3/2}$, the same as found by Kang et al.~(2020) when
  maximizing the finite-time growth of enstrophy. In addition, we also
  consider sharpness of an a priori estimate on the time evolution of
  $\| \mathbf{u}(t) \|_{L^4(\Omega)}$ by solving another PDE
  optimization problem and demonstrate that the upper bound in this
  estimate could be improved.
\end{abstract}



\section{Introduction}
\label{sec:intro}

This investigation concerns a systematic search for potentially
singular behavior in three-dimensional (3D) Navier-Stokes flows. By
formation of a ``singularity'' we mean the situation when an initially
smooth solution no longer satisfies the governing equation in the
classical (pointwise) sense. This so-called ``blow-up problem'' is one
of the key open questions in mathematical fluid mechanics and, in
fact, its importance for mathematics in general has been recognized by
the Clay Mathematics Institute as one of its ``millennium problems''
\cite{f00}.  Should such singular behavior indeed be possible in the
solutions of the 3D Navier-Stokes problem, it would invalidate this
system as a model of realistic fluid flows.  Questions concerning
global-in-time existence of smooth solutions remain open also for a
number of other flow models including the 3D Euler equations
\cite{gbk08} and some of the ``active scalar'' equations \cite{k10}.

We consider the incompressible Navier-Stokes system defined on the 3D
unit cube $\Omega = [0,1]^3$ with periodic boundary conditions
\begin{subequations}\label{eq:NS}
\begin{alignat}{2}
\partial_t\u + \u\cdot\bnabla\u + \bnabla p - \nu\laplacian\u & = 0 & &\qquad\mbox{in} \,\,\Omega\times(0,T], \\
\bnabla\cdot\u & = 0 & & \qquad\mbox{in} \,\,\Omega\times[0,T], \\
\u(0) & = \u_0, &   &
\end{alignat}
\end{subequations}
where the vector $\u = {[u_1, u_2, u_3]^T}$ is the velocity field, $p$
is the pressure, $\nu>0$ is the coefficient of kinematic viscosity and
$\u_0$ is the initial condition. The velocity gradient $\bnabla\u$ is
a tensor with components $[\bnabla\u]_{ij} = \partial_j u_i$,
$i,j=1,2,3$.  The fluid density $\rho$ is assumed constant and equal
to unity ($\rho=1$).

In our study an important role will be played by Lebesgue norms of the
velocity field
\begin{equation}
\|\u(t))\|_{L^q(\Omega)}  := \left( \int_\Omega |\u(t,\x)|^q \,d\x \right)^{\frac{1}{q}}, \qquad q \ge 1,
\label{eq:uLq}
\end{equation}
where ``$:=$'' means ``equal to by definition'', such that the kinetic
energy can be expressed as
\begin{equation}
\K(\u(t))  := \frac{1}{2} \|\u(t))\|_{L^2(\Omega)}.
\label{eq:K}
\end{equation}
Another important quantity is the enstrophy\footnote{{We note that
    unlike energy, cf.~\eqref{eq:K}, enstrophy is often defined
    without the factor of 1/2. However, for consistency with earlier
    studies belonging to this research program
    \cite{ap11a,ap13a,ap13b,ap16,Yun2018,KangYumProtas2020}, we choose
    to retain this factor here.}}
\begin{equation}
\E(\u(t))  :=  \frac{1}{2}\int_\Omega | \rot\u(t,\x) |^2 \,d\x
\label{eq:E}
\end{equation}
and the two quantities are related via the energy equation
\begin{equation}
\frac{d\K(\u(t))}{dt} = - \nu \E(\u(t)).
\label{eq:dKdt}
\end{equation}

While global in time existence of classical solutions of the
Navier-Stokes system \eqref{eq:NS} remains an open question, it is
known that suitably defined weak solutions, which need not satisfy the
Navier-Stokes system pointwise in space and time, but rather in a
certain integral sense only, exist globally in time \cite{l34}. An
important tool in the study of the global-in-time regularity of
classical (smooth) solutions are the so-called ``conditional
regularity results'' stating additional conditions which must be
satisfied by a weak solution in order for it to also be a smooth
solution, i.e., to satisfy the Navier-Stokes system in the classical
sense as well. One of the best known results of this type \cite{ft89}
is based on the enstrophy of the time-dependent velocity field $\u(t)$
and asserts that if the uniform bound
\begin{equation}
\label{eq:RegCrit_FoiasTemam}
\mathop{\sup}_{0 \leq t \leq T} \E(\u(t))  < \infty
\end{equation}
holds, then the regularity {and uniqueness of the solution $\u(t)$
  are} guaranteed up to time $T$ (to be precise, the solution remains
in a certain Gevrey class).

In the light of condition \eqref{eq:RegCrit_FoiasTemam} it is
important to characterize the largest growth of enstrophy possible in
Navier-Stokes flows.  Using \eqref{eq:NS} its rate of growth can be
expressed as $\frac{d\E(\u(t))}{dt} = -\nu\int_\Omega
|\laplacian\u|^2\,d\x + \int_{\Omega}
\u\cdot\nabla\u\cdot\laplacian\u\, d\x =: \R(\u(t))$ which is subject
to the following bound \cite{ld08,d09}
\begin{equation}
\frac{d\E}{dt} \leq \frac{27}{8\,\pi^4\,\nu^3} \E^3. 
\label{eq:dEdt_estimate_E}
\end{equation} 
By simply integrating the differential inequality in
\eqref{eq:dEdt_estimate_E} with respect to time we obtain the
finite-time bound
\begin{equation}
\E(\u(t)) \leq \frac{\E_0}{\sqrt{1 - \frac{27}{4\,\pi^4\,\nu^3}\,\E_0^2\, t}} 
\label{eq:Et_estimate_E0}
\end{equation}
which clearly becomes infinite at time $t_0 = 4\,\pi^4\,\nu^3 / (27\,
\E_0^2)$. Thus, based on estimate \eqref{eq:Et_estimate_E0}, it is not
possible to establish the boundedness of the enstrophy $\E(\u(t))$
required in condition \eqref{eq:RegCrit_FoiasTemam} and hence also the
regularity of solutions globally in time.

In addition to the enstrophy condition \eqref{eq:RegCrit_FoiasTemam},
another important conditional regularity result is given by the the
family of the Ladyzhenskaya-Prodi-Serrin conditions asserting that
Navier-Stokes flows $\u(t)$ are smooth and satisfy system
\eqref{eq:NS} in the classical sense provided that
\cite{KisLad57,Prodi1959,Serrin1962}
\begin{equation}
\u \in L^p([0,T];L^q(\Omega)), \quad 2/p+3/q = 1, \quad q > 3.
\label{eq:LPS}
\end{equation}
These conditions were recently generalized in \cite{Gibbon2018} to
include norms of the derivatives of the velocity field and to account
for velocity-pressure correlations in \cite{TranYuDritschel2021}. As
regards the limiting case with $q = 3$, the corresponding condition
was established in \cite{Escauriaza2003}
\begin{equation}
\u \in L^{\infty}([0,T];L^3(\Omega))
\label{eq:LPS3}
\end{equation}
and a related blow-up criterion was recently obtained in
\cite{Tao2020}.

Condition \eqref{eq:LPS} implies that should a singularity form in a
classical solution $\u(t)$ of the Navier-Stokes system \eqref{eq:NS}
at some finite time $0 < t_0 < \infty$, then necessarily
\begin{equation}
\lim_{t \rightarrow t_0} \int_0^t \| \u(\tau) \|_{L^q(\Omega)}^p \, d\tau \rightarrow \infty, \quad 2/p+3/q = 1, \quad q > 3.
\label{eq:LPSblowup}
\end{equation}
At the same time, the time evolution of the solution norm $\| \u(t)
\|_{L^q(\Omega)}$ on the time interval $[0,T]$ is subject to the some
a priori bounds valid also for Leray-Hopf weak solutions
\cite{Gibbon2018} which might involve singularities.  Such an estimate
was discussed in \cite{Constantin1991}
\begin{equation}
\int_0^T \| \u(\tau) \|_{L^q(\Omega)}^{\frac{4q}{3(q-2)}} \, d\tau \le  C \, \K_0^{\frac{2q}{3(q-2)}}, \qquad 2 \le q \le 6,
\label{eq:LPSbound}
\end{equation}
where $\K_0 := \K(\u_0)$ and $C>0$ is a generic constant whose
numerical value may vary between different estimates (since in
\cite{Constantin1991} this estimate is stated without an explicit
upper bound on the right-hand side (RHS), i.e., simply asserting the
finiteness of the expression on the left-hand side (LHS), estimate
\eqref{eq:LPSbound} is derived in Appendix \ref{sec:LPSbound}). We
note that the integral expressions in \eqref{eq:LPSblowup} and
\eqref{eq:LPSbound} differ in the exponent in the integrand which is
smaller by a factor of 2 in the latter case.  A related estimate,
known already to Leray \cite{l34}, concerns bounds on the rate of
growth of the $L^q$ norm and has the form
\cite{Giga1986,RobinsonSadowskiSilva2012,RobinsonSadowski2014}
\begin{equation}
\frac{1}{q} \frac{d}{dt} \| \u(t) \|_{L^q(\Omega)}^q \le C  \| \u(t) \|_{L^q(\Omega)}^{\frac{q(q-1)}{q-3}}, \qquad q > 3.
\label{eq:dLqdt}
\end{equation}
This estimate, which is analogous to \eqref{eq:dEdt_estimate_E}, makes
it possible to obtain lower bounds on $L^q$ norms of solutions
undergoing a hypothetical singularity formation in finite time.

We add that in the context of the inviscid Euler system a conditional
regularity result similar to \eqref{eq:RegCrit_FoiasTemam} and
\eqref{eq:LPS}--\eqref{eq:LPS3} is given by the Beale-Kato-Majda (BKM)
criterion \cite{bkm84}. Recently, finite-time singularity formation in
3D axisymmetric Euler flows on domains exterior to a boundary with
conical shape was proved in \cite{ElgindiJeong2018}.

While the blow-up problem is fundamentally a question in mathematical
analysis, a lot of computational studies have been carried out since
the mid-'80s in order to shed light on the hydrodynamic mechanisms
which might lead to singularity formation in finite time. Given that
such flows evolving near the edge of regularity involve formation of
very small flow structures, these computations typically require the
use of state-of-the-art computational resources available at a given
time. The computational studies focused on the possibility of
finite-time blow-up in the 3D Navier-Stokes and/or Euler system
include \cite{bmonmu83,ps90,b91,k93,p01,bk08,oc08,o08,ghdg08,
  gbk08,h09,opc12,bb12,opmc14,CampolinaMailybaev2018}, all of which
considered problems defined on domains periodic in all three
dimensions. The investigations \cite{dggkpv13,k13,gdgkpv14,k13b}
focused on the time evolution of vorticity moments and compared it
against bounds on these quantities obtained using rigorous analysis.
Recent computations \cite{Kerr2018} considered a ``trefoil''
configuration meant to be defined on an unbounded domain (although the
computational domain was always truncated to a finite periodic box). A
simplified semi-analytic model of vortex reconnection was recently
developed and analyzed based on the Biot-Savart law and asymptotic
techniques \cite{MoffattKimura2019a,MoffattKimura2019b}.  We also
mention the studies \cite{mbf08} and \cite{sc09}, along with
references found therein, in which various complexified forms of the
Euler equation were investigated. The idea of this approach is that,
since the solutions to complexified equations have singularities in
the complex plane, singularity formation in the real-valued problem is
manifested by the collapse of the complex-plane singularities onto the
real axis.  Overall, the outcome of these investigations is rather
inconclusive: while for the Navier-Stokes {system most of the} recent
computations do not offer support for finite-time blow-up, the
evidence appears split in the case of the Euler system.  In
particular, the studies \cite{bb12} and \cite{opc12} hinted at the
possibility of singularity formation in finite time. In this
connection we also highlight the {computational} investigations
\cite{lh14a,lh14b} in which blow-up was {documented} in axisymmetric
Euler flows on a bounded (tubular) domain. Recently, numerical
evidence for blow-up in solutions of the Navier-Stokes system in 3D
axisymmetric geometry with a degenerate variable diffusion coefficient
was provided in \cite{HouHuang2021}.

The related question of (non)uniqueness of solutions of the
Navier-Stokes system was considered in \cite{GuillodSverak2017} where
the authors focused on self-similar axisymmetric solutions
corresponding to initial data $\u_0$ with a singularity at the origin
chosen such that $\u_0$ is self-similar and does not belong to the
space $L^3(\RR^3)$. Nonunique solutions which do not satisfy
conditions \eqref{eq:LPS}--\eqref{eq:LPS3} were then found numerically
using the scale-invariance property to transform the Navier-Stokes
system to a nonlinear boundary-value problem. The problem of nonunique
solutions of 2D Euler equations corresponding to singular initial data
was recently tackled in \cite{BressanShen2021}.

A common feature of most of the aforementioned investigations was that
the initial data for the Navier-Stokes or Euler system was chosen in
an ad-hoc manner, based on some heuristic, albeit well-justified,
arguments. A new approach to the study of extreme, possibly singular,
behavior in fluid flows was ushered by Lu \& Doering who framed these
questions in terms of suitable variational optimization problems. In
\cite{l06,ld08} they showed that estimate \eqref{eq:dEdt_estimate_E}
is in fact sharp up to a numerical prefactor, in the sense that there
exists a family of velocity fields $\tuE \in H^2(\Omega)$
parameterized by their enstrophy $\E_0$ with the property that
$\frac{d}{dt}\E(\tuE) \sim \E_0^3$ as $\E_0 \rightarrow \infty$.
However, while these vector fields, which have the form of two
colliding vortex rings, saturate estimate \eqref{eq:dEdt_estimate_E}
{\em instantaneously}, the Navier-Stokes flows using these optimal
fields as the initial data feature rapid depletion of the rate of
enstrophy growth for $t > 0$ such that little enstrophy is produced
before it starts to decrease \cite{ap16} (for blow-up to occur,
enstrophy must be amplified at a sustained rate $\frac{d\E}{dt} \sim
\E^\alpha$, with $\alpha \in (2,3]$ for a sufficiently long time
\cite{KangYumProtas2020}). A research program where the sharpness of
various energy-type a priori estimates for one-dimensional (1D)
Burgers and two-dimensional (2D) Navier-Stokes flows was probed using
variational optimization formulations was pursued in
\cite{ap11a,ap13a,ap13b,ap16,Yun2018,ayala_doering_simon_2018}. While
these systems are known to be globally well-posed \cite{kl04},
questions about the sharpness of these estimates are quite pertinent
since these estimates are obtained in a similar way to the key
estimates \eqref{eq:Et_estimate_E0}, \eqref{eq:LPSblowup},
\eqref{eq:LPSbound} and \eqref{eq:dLqdt}.

The question whether enstrophy can become unbounded in finite time in
Navier-Stokes flows was investigated in \cite{KangYumProtas2020} by
finding optimal initial data $\tuET$ with fixed enstrophy $\E_0$ such
that the enstrophy is maximized at time $T$. This was done by solving
numerically a family of optimization problems 
\setcounter{problem}{-1}
\begin{problem}\label{pb:maxET}
  Given $\E_0, T \in\mathbb{R}_+$, find
\begin{align*}
\tuET & =  \mathop{\arg\max}_{\u_0 \in {\Q}_{\E_0}} \, \E(T), \quad \text{where} \\
{\Q}_{\E_0} & :=  \left\{\u\in H^1(\Omega)\,\colon\,\bnabla\cdot\u = 0, \; \E(\u) = \E_0 \right\},
\end{align*} 
\end{problem}
for a broad range of values of $\E_0$ and $T$. While no evidence was
found for unbounded growth of enstrophy in such extreme Navier-Stokes
flows, this study revealed the following approximate relation describing how the
largest attained enstrophy scales with the initial enstrophy in the
most extreme scenarios
\begin{equation}
\max_{T>0} \E(T) \approx 0.224 \, \E_0^{1.49}.
\label{eq:maxT_vs_E0}
\end{equation}
Interestingly, solution of an analogous maximization problem for 1D
viscous Burgers equation obtained in \cite{ap11a} produced extreme
flows which obey an essentially the same power-law relation as
\eqref{eq:maxT_vs_E0}, but with a different prefactor.

The goal of the present study is twofold: first, we will search for
initial data $\u_0$ which, subject to suitable constraints to be
defined below, might lead to unbounded growth of the integral in
\eqref{eq:LPSblowup} as $t \rightarrow t_0$, therefore signaling the
emergence of a singularity at time $t_0$; second, we will probe the
sharpness of the a priori estimate \eqref{eq:LPSbound} in terms of the
exponent of $\K_0$. More precisely, in regard to the second goal, the
objective is to verify whether the maximum of the quantity on the LHS
in \eqref{eq:LPSbound} achievable under the Navier-Stokes dynamics
\eqref{eq:NS} saturates the upper bound on the RHS, in the sense of
exhibiting the same scaling with the initial energy $\K_0$, which
would indicate that this estimate cannot be improved by reducing the
exponent of $\K_0$. To fix attention, we will consider these questions
for one only value of the parameter $q$. Concerning the first
question, we have found no evidence of unbounded growth required in
\eqref{eq:LPSblowup} to signal finite-time blow-up.  However,
interestingly, the families of the Navier-Stokes flows maximizing the
quantity $\int_0^T \| \u(\tau) \|_{L^4(\Omega)}^{8} \, d\tau$ for
different $T$ and different values of the constraint were found to
also follow a power-law relation with the same exponent as in
\eqref{eq:maxT_vs_E0} for the maximum growth of enstrophy. In regard
to the second question, we concluded that estimate \eqref{eq:LPSbound}
is not sharp, although the degree to which the upper bound
overestimates the growth of the expression on the LHS with $\K_0$ is
reduced as $T \rightarrow \infty$.

The structure of the paper is as follows: optimization problems
designed to probe the two questions mentioned above are stated in the
next section; then, in Section \ref{sec:approach} we introduce the
computational approach employed to solve these optimization problems;
our computational results are presented in Section \ref{sec:results},
whereas their discussion and conclusions are deferred to Section
\ref{sec:final}; some additional technical material is collected in
two appendices.

\section{Optimization Problem}
\label{sec:optim}

In this section we formulate optimization problems designed to provide
insights about the two questions stated in Introduction. For
concreteness, hereafter we will consider relations
\eqref{eq:LPSblowup} and \eqref{eq:LPSbound} with fixed values of the
indices $q=4$ and $p=8$. The reason for choosing these particular
values of $p$ and $q$ is our desire to work with integer-valued
indices, which will simplify numerical computations, while remaining
``close'' to the limiting critical case corresponding to $q=3$,
cf.~\eqref{eq:LPS3} (since this last condition is not given in terms
of an integral expression, it would need to be studied using methods
different from the approach developed here). We assume that with the
given initial data $\u_0$ the Navier-Stokes system \eqref{eq:NS}
admits classical solutions on the time interval $[0,T]$ and define the
quantities
\begin{subequations}\label{eq:PhiPsi}
\begin{align}
\Phi_T(\u_0) & := \frac{1}{T} \int_0^T \| \u(\tau) \|_{L^4(\Omega)}^8 \, d\tau, \label{eq:Phi} \\
\Psi_T(\u_0) & := \frac{1}{T} \int_0^T \| \u(\tau) \|_{L^4(\Omega)}^{8/3} \, d\tau, \label{eq:Psi} 
\end{align}
\end{subequations}
where $\u(t)$ is the solution of \eqref{eq:NS} with the initial
condition $\u_0$. These quantities correspond to the integral
expressions in \eqref{eq:LPSblowup} and \eqref{eq:LPSbound}, except
for the presence of the prefactor $T^{-1}$ whose role is to offset the
growth of the integrals which may occur for large $T$ even in the
absence of potentially singular events.

The idea for probing condition \eqref{eq:LPSblowup} is to formulate
and solve an optimization problem in order to find initial data $\u_0$
maximizing $\Phi_T(\u_0)$ for some $T>0$. However, for such an
optimization problem to be well defined, suitable constraints must be
imposed on $\u_0$ and a natural choice would be to require
$\|\u_0\|_{L^4(\Omega)} = B$ for some sufficiently large $0 < B <
\infty$ (the important question about the function space in which this
optimization problem should be posed is addressed below). Then, if a
hypothetical singularity is to occur at some time $t_0 > T$,
$\max_{\|\u_0\|_{L^4(\Omega) = B}} \Phi_T(\u_0)$ must become unbounded
as $T \rightarrow t_0$. Of course, a priori we do not know whether or
not a singularity may form, let alone at what time $t_0$, so condition
\eqref{eq:LPSblowup} can be probed by maximizing $\Phi_T(\u_0)$ for
increasing $T$ at a given value of $B$, and then repeating the process
for larger $B$. This approach is justified by upper bounds available
on the largest time $t_0$ when singularity might occur
\cite{Ohkitani2016}.

From the computational point of view, PDE-constrained optimization
problems are formulated most conveniently in a Hilbert space
\cite{pbh04}. While there exist solution approaches applicable in the
more general setting of Banach spaces, e.g., \cite{protas2008}, they
are significantly harder to use in practice. Given the form of our
constraint, we will therefore formulate the optimization problems in
the ``largest'' Sobolev space with Hilbert structure which is
contained in $L^4(\Omega)$. From the Sobolev embedding theorem in
dimension 3 \cite{af05}, we deduce
\begin{equation}
H^s(\Omega) \hookrightarrow L^4(\Omega), \qquad s \ge \frac{3}{4},
\label{eq:HsL4}
\end{equation}
where the Sobolev space $H^s(\Omega)$ is endowed with the norm $\| \z
\|_{H^{s}(\Omega)} = \| \z \|_{L^{2}(\Omega)} + \ell^{2s} \| \z
\|_{\dot{H}^{s}(\Omega)}$, $\forall \z \in H^{s}(\Omega)$, where $ \|
\z \|_{\dot{H}^{s}(\Omega)}= \| \Delta^{s/2} \z \|_{L^{2}(\Omega)}$ is
a semi-norm and $0 < \ell < \infty$ is the Sobolev parameter (while for
different values of $\ell$ the norms $\| \z \|_{H^{s}(\Omega)}$ are
equivalent, the choice of this parameter will play a role in numerical
computations, cf.~Section \ref{sec:results}). The fractional Laplacian
is defined in terms of the Fourier transform $\F$ as $ \Delta^{s/2} :=
\F^{-1} \left( | \k |^{s} \F \z \right)$, $s \in \RR$, where $\k \in
\ZZ^3$ is the wavevector. Thus, the largest Hilbert-Sobolev space
embedded in $L^4(\Omega)$ is the space $H^{3/4}(\Omega)$ and it will
provide the functional setting for our optimization problems.

We therefore arrive at the following 
\begin{problem}\label{pb:PhiL4}
  Given $B, T \in\mathbb{R}_+$ and the objective functional $\Phi_T(\u_0)$ from
  equation \eqref{eq:Phi}, find
\begin{align*}
\tuBT & =  \mathop{\arg\max}_{\u_0 \in {\L}_{B}} \, \Phi_T(\u_0), \quad \text{where} \\
 {\L}_{B} & :=  \left\{\u\in H^{3/4}(\Omega)\,\colon\,\bnabla\cdot\u = 0, \; \int_{\Omega} \u_0 \, d\x = \0,  \; \|\u_0\|_{L^4(\Omega)} = B \right\},
\end{align*} 
\end{problem}
\noindent
where the second condition in the definition of the constraint
manifold $ {\L}_B$ fixes the mean momentum since this quantity is
conserved under the evolution governed by the Navier-Stokes system
\eqref{eq:NS}.

Embedding \eqref{eq:HsL4} implies that $\forall \u \in
H^{3/4}(\Omega)$ $\|\u\|_{L^4(\Omega)} \le C \|\u\|_{H^{3/4}(\Omega)}$
and this allows us to re-express the constraint on the initial data
$\u_0$ in terms of its $H^{3/4}$ norm, which is quadratic in $\u_0$
and therefore easier to enforce in computations. This leads us to
\begin{problem}\label{pb:PhiH3/4}
  Given $S, T \in\mathbb{R}_+$ and the objective functional $\Phi_T(\u_0)$ from
  equation \eqref{eq:Phi}, find
\begin{align*}
\tuST & =  \mathop{\arg\max}_{\u_0 \in \H_{S}} \, \Phi_T(\u_0), \quad \text{where} \\
\H_{S} & :=  \left\{\u\in H^{3/4}(\Omega)\,\colon\,\bnabla\cdot\u = 0, \; \int_{\Omega} \u_0 \, d\x = \0,  \; \|\u_0\|_{\dot{H}^{3/4}(\Omega)} = S \right\}.
\end{align*} 
\end{problem}
We note that while Problem \ref{pb:PhiH3/4} is defined in the space
$H^{3/4}(\Omega)$, the last condition defining the constraint
manifold $\H_S$ is expressed in terms of the seminorm
$\|\cdot\|_{\dot{H}^{3/4}(\Omega)}$. This is done to ensure the
constraint manifold has a similar structure to the manifold $ {\L}_B$ in
Problem \ref{pb:PhiL4} and to the constraint manifold used in
\cite{KangYumProtas2020}.

A potential deficiency of Problem \ref{pb:PhiL4} is that the constraint
$\|\u_0\|_{L^4(\Omega)} = S > 0$ does not define a bounded set in the
space $H^{3/4}(\Omega)$, in the sense that one can construct a
sequence $\z_n \in H^{3/4}(\Omega)$, $n \in \ZZ$, such that $\forall
n$ $\|\z_n\|_{L^4(\Omega)} = S$ and $\lim_{n \rightarrow \infty}
\|\z_n\|_{H^{3/4}(\Omega)} = \infty$. However, while theoretically
possible, such behavior has not been observed in the computations
reported in Section \ref{sec:results}. 

As regards the second question we want to answer, concerning the
sharpness of estimate \eqref{eq:LPSbound}, given that the upper bound
in this estimate is expressed in terms of the initial energy $\K_0$, a
natural form of the corresponding optimization problem is given by
\begin{problem}\label{pb:PsiL2}
  Given $\K_0, T \in\mathbb{R}_+$ and the objective functional $\Psi_T(\u_0)$ from
  equation \eqref{eq:Psi}, find
\begin{align*}
\tuKT & =  \mathop{\arg\max}_{\u_0 \in \G_{\K_0}} \, \Psi_T(\u_0), \quad \text{where} \\
\G_{\K_0} & =  \left\{\u\in H^{3/4}(\Omega)\,\colon\,\bnabla\cdot\u = 0, \; \int_{\Omega} \u_0 \, d\x = \0,  \; \frac{1}{2}\|\u_0\|_{L^2(\Omega)} = \K_0 \right\},
\end{align*} 
\end{problem}
Our approach to solution of Problems \ref{pb:PhiL4}, \ref{pb:PhiH3/4}
and \ref{pb:PsiL2} is described next.

\section{Computational Approach}
\label{sec:approach}

In this section we describe our approach to solution of optimization
problems \ref{pb:PhiL4}, \ref{pb:PhiH3/4} and \ref{pb:PsiL2} for given
values of $B$, $S$ or $\K_0$ and $T$. We adopt an
``optimize-then-discretize'' approach \cite{g03} in which a gradient
method is first formulated in {the} infinite-dimensional (continuous)
setting and only then the resulting {equations and} expressions are
discretized for the purpose of numerical solution. A similar approach
was recently used to solve the problem of determining the maximum
growth of enstrophy in \cite{KangYumProtas2020} with the corresponding
1D problem addressed earlier in \cite{ap11a}. To make the present
paper self-contained, we recall key elements of the solution approach
from \cite{KangYumProtas2020}. However, there are also some important
differences resulting from the functional setting and the nature of
the constraints in Problems \ref{pb:PhiL4}, \ref{pb:PhiH3/4} and
\ref{pb:PsiL2} which we highlight. We also mention the Riemannian
aspects of the optimization problems \cite{ams08}. In our presentation
below we first focus on solving Problem \ref{pb:PhiH3/4} as it
arguably has the simplest structure and then discuss the modifications
required to solve Problems \ref{pb:PhiL4} and \ref{pb:PsiL2}.

\subsection{Discrete Gradient Flow}
\label{sec:gradflow}

Problem \ref{pb:PhiH3/4} is Riemannian in the sense that the maximizer
$\tuST$ must be contained on a constraint manifold $\H_S$
\cite{ams08}. In order to locally characterize this manifold, we
construct the tangent space $\T_{\z}\H_S$ at some point $\z \in \H_S$.
The fixed-norm constraints can be expressed in terms of the function
$F_X \; : \; H^{3/4}(\Omega) \rightarrow \RR_+$, $F_X := \|\z\|_X$,
where $X = L^4(\Omega), \dot{H}^{3/4}(\Omega), L^2(\Omega)$
respectively in Problems \ref{pb:PhiL4}, \ref{pb:PhiH3/4} and
\ref{pb:PsiL2}. Then, the subspace tangent to the manifold defined in
the space $H^{3/4}(\Omega)$ by the relation $F_X(\z) = S$ is given by
the condition $\forall \z' \in H^{3/4}(\Omega)$ $F'(\z;\z') = \langle
\bnabla F_X(\z),\z' \rangle_{H^{3/4}(\Omega)} = 0$ which also defines
the element $\bnabla F_X(\z)$ orthogonal to the subspace.  Thus, since
in Problem \ref{pb:PhiH3/4} we have $F_X(\z) = F_{\dot{H}^{3/4}}(\z)
:= \| \z \|_{\dot{H}^{3/4}}$, the tangent space to the manifold $\H_S$
is defined as
\begin{align}
\T_{\z}\H_S & := \left\{ \v \in H^{3/4}(\Omega) \, : \, \bnabla\cdot\v = 0, \; \int_{\Omega} \v \, d\x = \0,  \; \left\langle \bnabla F_{\dot{H}^{3/4}}(\z), \v \right\rangle_{H^{3/4}(\Omega)} = 0 \right\},
\label{eq:TzHS} \\
& \text{where} \qquad  \left\langle \bnabla F_{\dot{H}^{3/4}}(\z), \z' \right\rangle_{H^{3/4}(\Omega)} = 
\langle \z, \z' \rangle_{\dot{H}^{3/4}(\Omega)}, \quad \forall \z' \in H^{3/4}(\Omega) \nonumber 
\end{align}
(we note that in general $\bnabla F_{\dot{H}^{3/4}}(\z) \neq \z$ since
the constraint is defined in terms of the semi-norm).  

The maximizer $\tuST$ can then be found as $\tuST = \lim_{n\rightarrow
  \infty} \uST^{(n)}$ using the following iterative procedure
representing a discretization of a gradient flow projected on $\H_S$
\begin{equation}
\begin{aligned}
\uST^{(n+1)} & =  \RR_{\H_S}\left(\;\uST^{(n)} + \tau_n \, \PP_{\T_n}\nabla\Phi_T\left(\uST^{(n)}\right)\;\right), \\ 
\uST^{(1)} & =  \u^0.
\end{aligned}
\label{eq:desc}
\end{equation}
Here $\uST^{(n)}$ is an approximation of the maximizer obtained at the
$n$-th iteration, $\u^0$ is the initial guess, $\PP_{\T_n} \; : \;
H^{3/4}(\Omega) \rightarrow \T_n := \T_{\uST^{(n)}}\H_S$ is an
operator representing projection onto the tangent subspace
\eqref{eq:TzHS} at the $n$th iteration, $\tau_n$ is the length of the
step whereas $\RR_{\H_S} \: : \; \T_n \rightarrow \H_S$ is a
retraction from the tangent space to the constraint manifold
\cite{ams08}.  A key element of the iterative procedure
\eqref{eq:desc} is the evaluation of the gradient $\nabla\Phi_T$ of
the objective functional $\Phi_T$, cf.~\eqref{eq:Phi}, representing
its (infinite-dimensional) sensitivity to perturbations of the initial
data $\u_0$ {in the governing system \eqref{eq:NS}. We emphasize that
  it is essential for the gradient to possess} the required
regularity, namely, $\nabla\Phi_T(\u_0) \in H^{3/4}(\Omega)$.

The first step to determine the gradient $\nabla\Phi_T$ is to consider
the G\^{a}teaux (directional) differential $\Phi'_T(\u_0;\cdot) \; : \;
H^{3/4}(\Omega) \rightarrow \RR$ of the objective functional $\Phi_T$ defined
as $\Phi'_T(\u_0;\u_0') := \lim_{\epsilon \rightarrow 0}
\epsilon^{-1}\left[\Phi_T(\u_0+\epsilon \u_0') - \Phi_T\right]$ for some
arbitrary perturbation $\u_0' \in H^{3/4}(\Omega)$. The gradient
$\nabla\Phi_T$ can {then} be extracted from the G\^{a}teaux
differential $\Phi'_T(\u_0;\u_0')$ recognizing that, when viewed as a
function of its second argument, this differential is a bounded linear
functional on the space $H^{3/4}(\Omega)$ and we can therefore invoke the
Riesz representation theorem \cite{l69}
\begin{equation}
\Phi'_T(\u_0;\u_0')
= \Big\langle \nabla^{L_2}\Phi_T, \u_0' \Big\rangle_{L^2(\Omega)} = \Big\langle \nabla\Phi_T, \u_0' \Big\rangle_{H^{3/4}(\Omega)},
\label{eq:riesz}
\end{equation}
where the gradient $\nabla\Phi_T$ is the Riesz representer in the
function space $H^{3/4}(\Omega)$. In \eqref{eq:riesz} we also formally
defined the gradient $\nabla^{L_2}\Phi_T$ determined with respect to
the $L^2$ topology as it will be useful in subsequent computations.
Given the definition of the objective functional in \eqref{eq:Phi},
its G\^{a}teaux differential can be expressed as
\begin{equation}
\Phi'_T(\u_0;\u_0') = \frac{8}{T} \int_0^T \left( \|\u(t)\|_{L^4(\Omega)}^4 \, 
\int_{\Omega} |\u(t,\x)|^2 \u(t,\x) \cdot \u'(t,\x) \, d\x \right) \, dt,
\label{eq:dPhiT}
\end{equation}
where the perturbation field $\u' = \u'(t,\x)$ is a solution of the
Navier-Stokes system linearized around the trajectory corresponding to
the initial data $\u_0$ \cite{g03}, i.e.,
\begin{subequations}
\label{eq:lNSE3D}
\begin{align}
 \L\begin{bmatrix} \u' \\ p' \end{bmatrix} := 
& \begin{bmatrix}
\partial_{t}\u'+\u'\cdot\bnabla\u+\u\cdot\bnabla\u'+\bnabla p'-\nu\bDelta\u' \\
\bnabla\cdot\u'
\end{bmatrix} = \begin{bmatrix} \mathbf{0} \\ 0\end{bmatrix}, \label{eq:lNSE3Da} \\
 \u'(0)= &\u_0' \label{eq:lNSE3Db}
\end{align}
\end{subequations}
which is subject to the periodic boundary conditions and where $p'$ is
the perturbation of the pressure.

We note that expression \eqref{eq:dPhiT} for the G\^{a}teaux
differential is not yet consistent with the Riesz form
\eqref{eq:riesz}, because the perturbation $\u_0'$ of the initial data
does not appear in it explicitly as a factor, but is instead hidden as
the initial {condition} in the linearized problem,
cf.~\eqref{eq:lNSE3Db}. In order to transform \eqref{eq:dPhiT} to the
Riesz form, we introduce the {\em adjoint states} $\u^* \; : \;
[0,T]\times\Omega \rightarrow \RR^3$ and $p^* \; : \;
[0,T]\times\Omega \rightarrow \RR$, and the following duality-pairing
relation
\begin{equation}
\begin{aligned}
\left( \L\begin{bmatrix} \u' \\ p' \end{bmatrix}, \begin{bmatrix} \u^* \\ p^* \end{bmatrix} \right)
:= & \int_0^T \int_{\Omega} \L\begin{bmatrix} \u' \\ p' \end{bmatrix} \cdot \begin{bmatrix} \u^* \\ p^* \end{bmatrix} \, d\x \, dt 
= \overbrace{\left( \begin{bmatrix} \u' \\ p' \end{bmatrix}, \L^*\begin{bmatrix} \u^* \\ p^* \end{bmatrix}\right)}^{\Phi'_T(\u_0;\u_0')} + \\
\phantom{=} & \int_\Omega \u'(T,\x)\cdot\u^*(T,\x)  \,d\x - 
\int_\Omega \u'(0,\x)\cdot\u^*(0,\x)  \,d\x = 0,
\end{aligned}
\label{eq:dual}
\end{equation}
where ``$\cdot$'' in the first integrand expression denotes the
Euclidean dot product evaluated at $(t,\x)$. {Performing}
integration by parts with respect to {both space and time then}
allows us to define the {\em adjoint system} as
\begin{subequations}
\label{eq:aNSE3D}
\begin{align}
 \L^*\begin{bmatrix} \u^* \\ p^* \end{bmatrix} := 
& \begin{bmatrix}
-\partial_{t}\u^*-\left[\bnabla\u^*+\left(\bnabla\u^{*}\right)^T\right]\u-\bnabla p^*-\nu\bDelta\u^* \\
-\bnabla\cdot\u^*
\end{bmatrix}  = \begin{bmatrix} \f  \\ 0\end{bmatrix}, \label{eq:aNSE3Da} \\
\text{where} \quad \f(t,\x)  := & \frac{8}{T}  \|\u(t)\|_{L^4(\Omega)}^4 \,  |\u(t,\x)|^2 \u(t,\x), \qquad \x\in \Omega, \quad t \in [0,T], \label{eq:aNSE3Df} \\
 \u^*(T)= & \0  \label{eq:aNSE3Db}
\end{align}
\end{subequations}
which is also subject to the periodic boundary conditions. We note
that in identity \eqref{eq:dual} all boundary terms resulting from
integration by parts with respect to the space variable vanish due to
the periodic boundary conditions. The term $\int_\Omega
\u'(T,\x)\cdot\u^*(T,\x) \,d\x$ {resulting from integration by parts
  with respect to time vanishes because of the homogeneous terminal
  condition \eqref{eq:aNSE3Db} such that with the judicious} choice of
the source term \eqref{eq:aNSE3Df} identity \eqref{eq:dual} implies
\begin{equation}
\Phi'_T(\u_0;\u_0') = \int_\Omega \u'_0(\x)\cdot\u^*(0,\x)  \,d\x.
\label{eq:dET2}
\end{equation}
Applying the first equality in Riesz relations \eqref{eq:riesz} to
\eqref{eq:dET2} we obtain the $L^2$ gradient as
\begin{equation}
\nabla^{L^2}\Phi_T = \u^*(0).
\label{eq:gradL2}
\end{equation}

Our Sobolev gradient $\nabla\Phi_T(\u_0)$ is defined in a fractional
Sobolev space $H^{3/4}(\Omega)$. However, since system \eqref{eq:NS}
is defined on a periodic domain $\Omega$, such a gradient can be
determined in a similar manner to the case of a Sobolev space with an
integer differentiability index \cite{pbh04}. We thus proceed by
identifying the G\^{a}teaux differential in \eqref{eq:dET2} with the
$H^{3/4}$ inner product. {Then, recognizing that
  the perturbations $\u_0'$ are arbitrary, we obtain the following
  fractional elliptic boundary-value problem
\begin{equation}
\left[ \Id \, - \,\ell^{3/2} \,\bDelta^{3/4} \right] \nabla\Phi_T(\u_0)
= \nabla^{L_2} \Phi_T(\u_0)   \qquad \text{in} \ \Omega 
\label{eq:gradH3/4}
\end{equation}
subject to the periodic boundary conditions, which must be solved to
determine $\nabla\Phi_T$.  System \eqref{eq:gradH3/4} is conveniently
solved in the Fourier space where it takes the form
\begin{subequations}
\label{eq:FgradH3/4}
\begin{align}
\left[ 1 + \,\ell^{3/2} \,|\k|^{3/2} \right] \left[\widehat{\nabla\Phi_T(\u_0)}\right]_\k 
& = \left[\widehat{\nabla^{L^2}\Phi_T(\u_0)}\right]_\k,    \qquad \k \in \ZZ^3 \setminus \0, \label{eq:FgradH3/4k} \\
\left[\widehat{\nabla\Phi_T(\u_0)}\right]_\0 & = \0,  \label{eq:FgradH3/40}
\end{align}
\end{subequations}
in which $\left[\widehat{\z}\right]_\k \in \CC^3$ denotes the Fourier
coefficient of the vector field $\z$ corresponding to the wavevector
$\k$. We remark that \eqref{eq:FgradH3/40} ensures that the Sobolev
gradient $\nabla\Phi_T(\u_0)$ satisfies the zero-mean condition in
Problem \ref{pb:PhiH3/4} (including this condition in system
\eqref{eq:FgradH3/4} is equivalent to projecting the resulting
gradient on the subspace defined by this condition).

The gradient fields $\nabla^{L_2}\Phi_T$ and $\nabla\Phi_T$ can be
interpreted as infinite-dimensional sensitivities of the objective
functional $\Phi_T$, cf.~\eqref{eq:Phi}, with respect to perturbations
of the initial data $\u_0$. While these two gradients point towards
the same local maximizer, they represent distinct ``directions'',
since they are defined with respect to different {norms} ($L^2$
vs.~$H^{3/4}$). As shown by \cite{pbh04}, extraction of gradients in
spaces of smoother functions such as $H^{3/4}(\Omega)$ can be
interpreted as low-pass filtering of the $L^2$ gradients with the
parameter $\ell_1$ acting as the cut-off length-scale. {Although
  Sobolev gradients obtained with different $0 < \ell_1 < \infty$ are
  equivalent, in the precise sense of norm equivalence \cite{b77}, in
  practice the value of $\ell_1$ tends to have a significant effect on
  the rate of convergence of gradient iterations \eqref{eq:desc}
  \cite{pbh04}} and the choice of its numerical value will be
discussed in Section \ref{sec:numer}. We emphasize that, while the
$H^{3/4}$ gradient is used exclusively in the actual computations,
cf.~\eqref{eq:desc}, the $L^2$ gradient is computed first as an
intermediate step.

Evaluation of the $L^2$ gradient at a given iteration {via
  \eqref{eq:gradL2}} requires solution of the Navier-Stokes system
\eqref{eq:NS} followed by solution of the adjoint system
\eqref{eq:aNSE3D}. We note that this system is a linear problem with
coefficients and the source term determined by the solution of the
Navier-Stokes system obtained earlier during the iteration. The
adjoint system \eqref{eq:aNSE3D} is a {\em terminal} value problem,
implying that it must be integrated {\em backwards} in time from $t=T$
to $t=0$ (since the term with the time derivative has a negative sign,
this problem is well posed).  Once the $L^2$ gradient is determined
using \eqref{eq:gradL2}, the corresponding Sobolev $H^{3/4}$ gradient
can be obtained by solving problem \eqref{eq:gradH3/4}. We add that
the thus computed Sobolev gradient satisfies the divergence-free
condition by construction, i.e., $\bnabla\cdot(\nabla\Phi_T) = 0$.

\subsection{Projection, Retraction and Arc-Maximization}
\label{eq:proj}

The projection operator $\PP_{\T_n}$ appearing in \eqref{eq:desc} is
defined as \cite{l69}, cf.~\eqref{eq:TzHS},
\begin{equation}
\forall_{\z \in H^{3/4}(\Omega)} \qquad \PP_{\T_n} \z := 
\z - \frac{\left\langle \z, \bnabla F_{\dot{H}^{3/4}}\left(\uST^{(n)}\right) \right\rangle_{H^{3/4}(\Omega)}}{\left\| \bnabla F_{\dot{H}^{3/4}}\left(\uST^{(n)}\right) \right\|_{H^{3/4}(\Omega)}}\, \bnabla F_{\dot{H}^{3/4}}\left( \uST^{(n)}\right).
\label{eq:P}
\end{equation}
As can be readily verified, it preserves both the divergence-free and
zero-mean conditions. The projection defined in \eqref{eq:P} can be
applied with obvious modifications consisting in changes of the norm
and the inner product to Problem \ref{pb:PsiL2}, but not to Problem
\ref{pb:PhiL4}. Expression for the projection operator in Problem
\ref{pb:PhiL4} will be discussed in Section \ref{sec:projL4}.

The retraction operator is defined as the normalization \cite{ams08}
\begin{equation}
\forall_{\z \in \T_n} \qquad  \RR_{\H_S}(\z) := \frac{S}{\|\z\|_{\dot{H}^{3/4}(\Omega)}}\,\z
\label{eq:R}
\end{equation}
which clearly also preserves the divergence-free and zero-mean
properties of the argument. The retraction defined in \eqref{eq:R} can
be applied with obvious adjustments to Problems \ref{pb:PhiL4} and
\ref{pb:PsiL2}. Projection of the gradient $\nabla\Phi_T(\u_0)$ onto
the tangent subspace $\T_n$ via \eqref{eq:P} followed by retraction
\eqref{eq:R} to the constraint manifold $\H_S$ are illustrated
schematically in Figure \ref{fig:HS}.
\begin{figure}
\begin{center}
\includegraphics[width=0.6\textwidth]{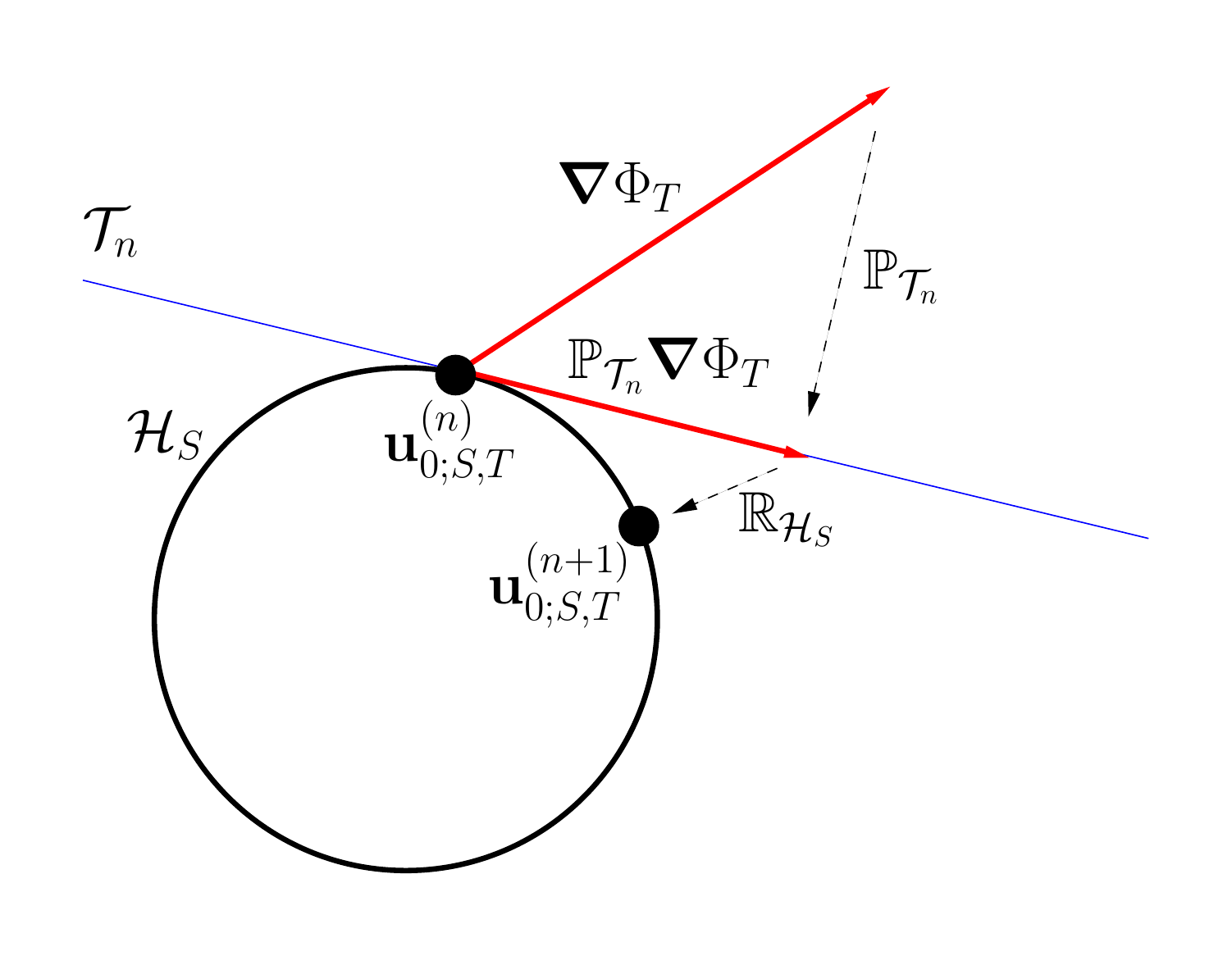}
\caption{Schematic representation of projection \eqref{eq:P} applied
  to the gradient $\bnabla\Phi_T(\u_0)$ followed by retraction
  \eqref{eq:R} to the manifold $\H_S$. }
\label{fig:HS}
\end{center}
\end{figure}

The step size $\tau_n$ in algorithm \eqref{eq:desc} is computed by
solving the problem
\begin{equation}\label{eq:tau_n}
\tau_n = \mathop{\argmax}_{\tau>0} \Phi_T\left( \RR_{\H_S}\left(\;\uST^{(n)} + \tau \, \PP_{\T_n}\nabla\Phi_T\left(\uST^{(n)}\right)\;\right) \right)
\end{equation}
which is done using a suitable derivative-free {approach, such as a
  variant of Brent's algorithm \cite{nw00,numRecipes}}.  Equation
\eqref{eq:tau_n} can be interpreted as a modification of {the standard
  line-search problem} where maximization is performed following an
arc (a geodesic in the limit of infinitesimal step sizes) lying on the
constraint manifold $\H_S$, rather than {along} a straight line.

\subsection{Projection on Tangent Subspace in Problem \ref{pb:PhiL4}}
\label{sec:projL4}

In Problem \ref{pb:PhiL4} the constraint is defined in terms of the
function $F_{L^4}(\z) := \| \z \|_{L^4(\Omega)}$, such that the
subspace tangent to the manifold $ {\L}_B$ is given by the condition
$\langle \bnabla F_{L^4}(\z), \z' \rangle_{H^{3/4}(\Omega)} = 0$,
$\forall \z' \in H^{3/4}(\Omega)$, where $\left\langle \bnabla
  F_{L^4}(\z), \z' \right\rangle_{H^{3/4}(\Omega)} = \langle |\z|^2\z,
\z' \rangle_{\dot{H}^{3/4}(\Omega)}$. We note that given the
nonlinearity of the term $|\z|^2\z$, the element $\bnabla F_{L^4}(\z)$
does not in general satisfy the divergence-free and zero-mean
conditions, even if they are satisfied by $\z$. Thus, projection
\eqref{eq:P} must be modified such that the result is both
divergence-free and has zero mean which is done as follows
\begin{align}
\forall \, {\z \in H^{3/4}(\Omega)} \qquad \PP_{\T_n} \z & := 
\z - \frac{\left\langle \z, \bnabla F_{L^4}\left(\uBT^{(n)}\right) \right\rangle_{H^{3/4}(\Omega)}}{\left\langle \overline{\bnabla F_{L^4}\left( \uBT^{(n)}\right)}, \bnabla F_{L^4}\left(\uBT^{(n)}\right) \right\rangle_{H^{3/4}(\Omega)}}\, \overline{\bnabla F_{L^4}\left( \uBT^{(n)}\right)}
\label{eq:PL4} \\
& \text{where} \quad  \overline{\v} := \v - \bnabla \Delta^{-1} (\bnabla \cdot \v) - \int_{\Omega} \v \, d\x. \nonumber
\end{align}

\subsection{Numerical Implementation}
\label{sec:numer}

The approach described in Sections
\ref{sec:gradflow}--\ref{sec:projL4} is implemented as described in
detail in \cite{KangYumProtas2020}. Here we summarize key elements of
the numerical methodology and refer the reader to
\cite{KangYumProtas2020} for further particulars.  Evaluation of the
objective functionals \eqref{eq:Phi}--\eqref{eq:Psi} requires solution
of the Navier-Stokes system \eqref{eq:NS} on the time interval $[0,T]$
with the given initial data $\u_0$, whereas determination of the $L^2$
gradient \eqref{eq:gradL2} requires solution of the adjoint system
\eqref{eq:aNSE3D}.  These two PDE systems are solved numerically with
an approach combining a pseudo-spectral approximation of spatial
derivatives with a {fourth-order} semi-implicit Runge-Kutta method
\cite{NumRenaissance} used to discretize these problems in time.  In
the evaluation of the nonlinear term in \eqref{eq:NS} and the terms
with non-constant coefficients in \eqref{eq:aNSE3D} dealiasing is
performed using the Gaussian filtering approach proposed in
\cite{hl07}.  The velocity field $\u = \u(t,\x)$ needed to evaluate
the coefficients and the source term in the adjoint system
\eqref{eq:aNSE3D} is saved at discrete time levels during solution of
the Navier-Stokes system \eqref{eq:NS}. In the definition of the
Sobolev gradient in \eqref{eq:gradH3/4}--\eqref{eq:FgradH3/4} we set
$\ell = 2$ which was found by trial and error to maximize the rate of
convergence of iterations \eqref{eq:desc}.  Massively parallel
implementation based on MPI and using the {\tt fftw} routines
\cite{fftw} to perform Fourier transforms allowed us to {employ}
resolutions varying from $128^3$ to $512^3$ in cases with low and high
values of the constraints, respectively. In the latter cases solution
of Problems \ref{pb:PhiL4}, \ref{pb:PhiH3/4} and \ref{pb:PsiL2} for an
intermediate length $T$ of the time interval typically required a
computational time of $\O(10^2)$ hours on $\O(10^2)$ CPU cores.  The
computational results presented in the next section have been
thoroughly validated using strategies described in
\cite{KangYumProtas2020} to ensure they are converged with respect to
refinement of the different numerical parameters.

Problems \ref{pb:PhiL4}, \ref{pb:PhiH3/4} and \ref{pb:PsiL2} are
non-convex and as such may admit multiple local maximizers. With the
gradient-based approach \eqref{eq:desc}, which relies on local
common
information only, we cannot assert whether the maxima we find are
global or not. In order to find as many local maxima as possible, for
each set of parameters $T$ and $B$, $S$ or $\K_0$ we solve Problems
\ref{pb:PhiL4}, \ref{pb:PhiH3/4} and \ref{pb:PsiL2} using different
initial guesses $\u^0$. For example, for Problem \ref{pb:PhiL4} we fix
the value of the constraint $B$ and then the corresponding branch of
maximizing solutions is obtained by solving the problem for a sequence
of (increasing or decreasing) values of $T$ using the optimal solution
$\tuBT$ obtained for the previous value of $T$ as the initial guess
$\u^0$. Then, another branch of maximizing solutions is obtained by
repeating this process for a different value of the constraint $B$.
We refer the reader to \cite{KangYumProtas2020} for further details of
this ``continuation'' approach. In addition, to make this search more
exhausting, we have also used various random initial guesses and
the optimal initial conditions found in \cite{KangYumProtas2020} as
the initial guess $\u^0$.

\section{Results}
\label{sec:results}

In this section we first discuss the results obtained by solving
Problems \ref{pb:PhiL4} and \ref{pb:PhiH3/4} designed to search for
initial data $\u_0$ that would trigger the appearance of a singularity
in finite time. Next, we present the results obtained by solving
Problem \ref{pb:PsiL2} defined to probe the sharpness of estimate
\eqref{eq:LPSbound}. In these calculations we set $\nu=0.01$ which is
the same value as used in earlier studies of closely-related problems
\cite{ld08,ap16,KangYumProtas2020}. In addition to other diagnostic
quantities, in our analysis of the different flows we will also
consider their componentwise enstrophies $\E_i(\u(t))$, $i=1,2,3$,
associated with the three coordinate directions and defined as
\begin{equation}
\E_i(\u(t)) := \int_{\Omega} \left| \left(\bnabla \times \u(t) \right) \cdot \e_i \right|^2\, d\x, \quad i=1,2,3,
\label{eq:Ei}
\end{equation}
where $\e_1$, $\e_2$, $\e_3$ are the unit vectors of the Cartesian
coordinate system and we have the obvious identity $\forall t \ \
\E(\u(t)) = \sum_{i=1}^3 \E_i(\u(t))$.

\subsection{Flows Obtained as Solutions of Problems  \ref{pb:PhiL4} and \ref{pb:PhiH3/4}}
\label{sec:resultsP1P2}

Solution of Problems \ref{pb:PhiL4} and \ref{pb:PhiH3/4} has yielded
two distinct maximizing branches in each case and representative
solutions are shown in terms of the time evolution of the norm
$\|\u(t)\|_{L^4(\Omega)}^4$ in Figures \ref{fig:L4vst}a and
\ref{fig:L4vst}b, respectively. Both figures show evolutions obtained
with the largest considered values of the constraints $B$ and $S$ for
``short'' and ``long'' optimization windows $T$. In regard to Problem
\ref{pb:PhiL4}, we see that for solutions from both maximizing
branches the quantity $\|\u(t)\|_{L^4(\Omega)}^4$ exhibits a
significant transient growth with larger maximum values $\max_{0 \le t
  \le T} \|\u(t)\|_{L^4(\Omega)}^4$ achieved for shorter optimization
windows $T$. On the other hand, for Problem \ref{pb:PhiH3/4} we note
that the norm $\|\u(t)\|_{L^4(\Omega)}^4$ exhibits monotone decrease
with time for maximizing solutions from both branches, cf.~Figure
\ref{fig:L4vst}b. These flows are quite similar to each other in terms
of the evolution of the norm $\|\u(t)\|_{L^4(\Omega)}^4$ and,
moreover, show weak dependence on the length $T$ of the optimization
window (in the sense that the flows obtained by solving Problem
\ref{pb:PhiH3/4} with $T_1$ and $T_2$ such that $T_1 < T_2$ exhibit a
similar evolution of $\|\u(t)\|_{L^4(\Omega)}^4$ for $t \in [0,T_1]$).

\begin{figure}
\begin{center}
\mbox{\subfigure[]{\includegraphics[width=0.5\textwidth]{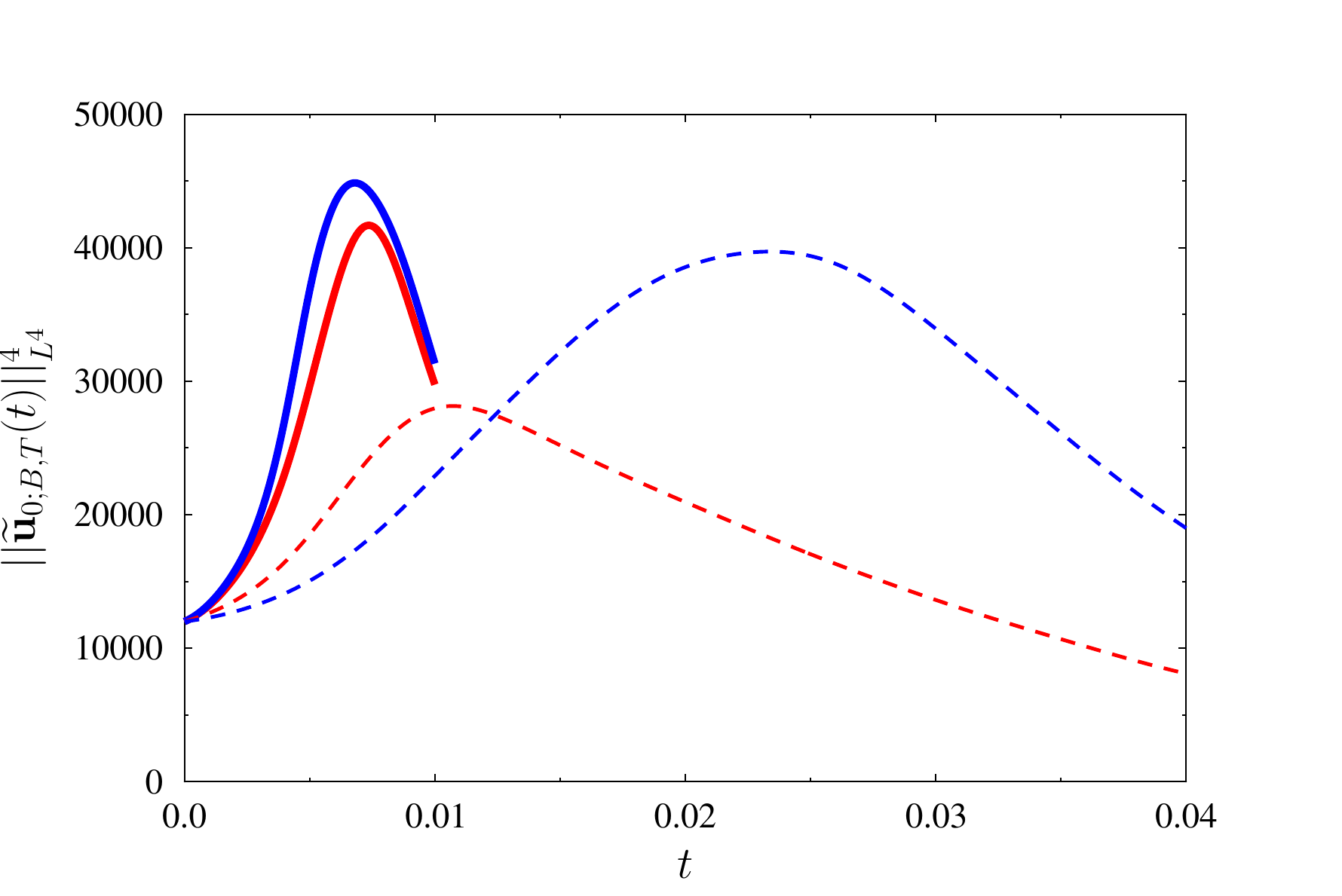}}\quad
\subfigure[]{\includegraphics[width=0.5\textwidth]{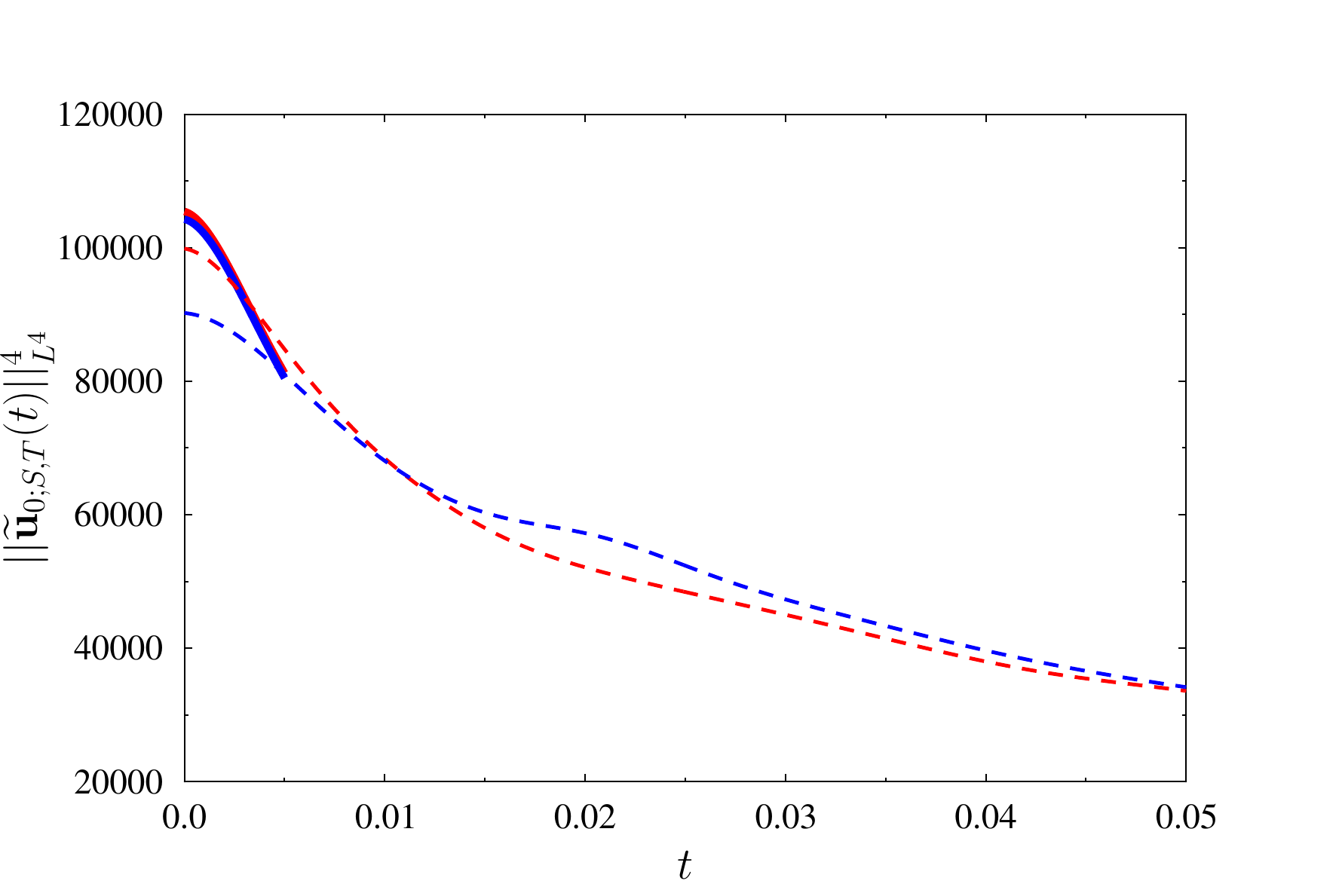}}}
\caption{Time evolution of $\| \tu (t) \|_{L^4}^{4}$ in the
  Navier-Stokes flows with optimal initial conditions obtained by
  solving (a) Problem \ref{pb:PhiL4} with $B^4=12{,}000$ and (b)
  Problem \ref{pb:PhiH3/4} with $S^2=2{,}000$ using two different
  lengths $T$ of the optimization window (in each case the results are
  shown on the time interval $[0,T]$ where optimization was
  performed). In panel (a) the blue and red lines correspond to the
  partially symmetric and asymmetric branches, whereas in panel (b)
  they correspond to the symmetric and two-component branches.}
\label{fig:L4vst}
\end{center}
\end{figure}

The maximizing branches obtained by solving Problems \ref{pb:PhiL4}
and \ref{pb:PhiH3/4} with five different values of the constraints $B$
and $S$ are shown in terms of the dependence of the maximum values of
the objective functional \eqref{eq:Phi} on the length $T$ of the
optimization window in Figures \ref{fig:maxLPSvsT}a and
\ref{fig:maxLPSvsT}b, respectively. The presence of two distinct
branches for each value of the constraint $B$ and $S$ is clearly
evident, although the differences are small for solutions of Problem
\ref{pb:PhiH3/4}, cf.~Figure \ref{fig:maxLPSvsT}b. We note that as
regards solutions of Problem \ref{pb:PhiL4}, for each value of the
constraint $B$, the largest values of the objective functional
$\Phi_T(\tuBT)$ are for both branches attained on optimization windows
with length $T$ decreasing with $B$, cf.~Figure \ref{fig:maxLPSvsT}a.
On the other hand, for solutions of Problem \ref{pb:PhiH3/4} obtained
with a fixed value of the constraint $S$, the maxima of the objective
functional $\Phi_T(\tuST)$ are in all cases decreasing functions of
the length $T$ of the optimization window.

\begin{figure}
\begin{center}
\mbox{\subfigure[]{\includegraphics[width=0.5\textwidth]{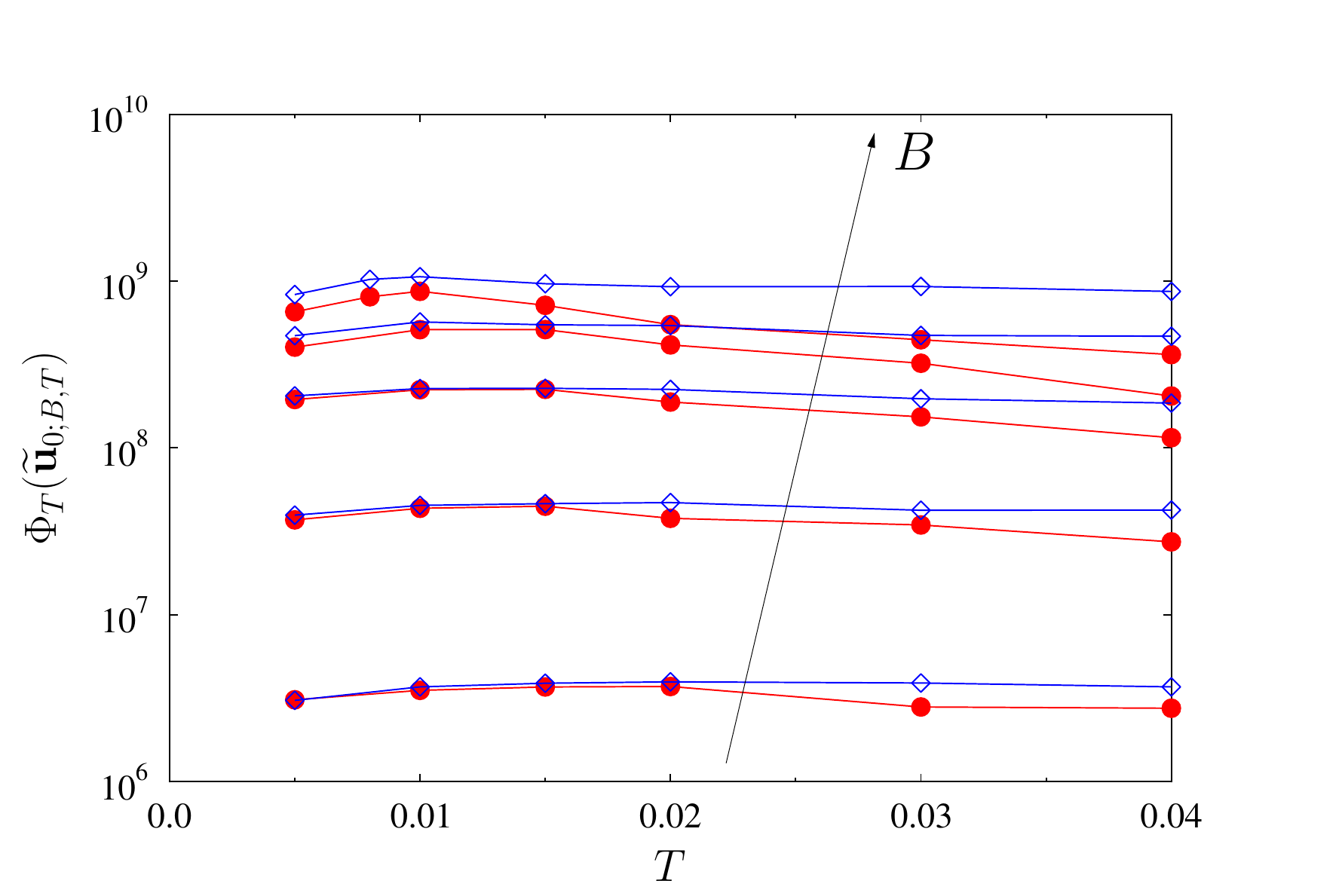}}\quad
\subfigure[]{\includegraphics[width=0.5\textwidth]{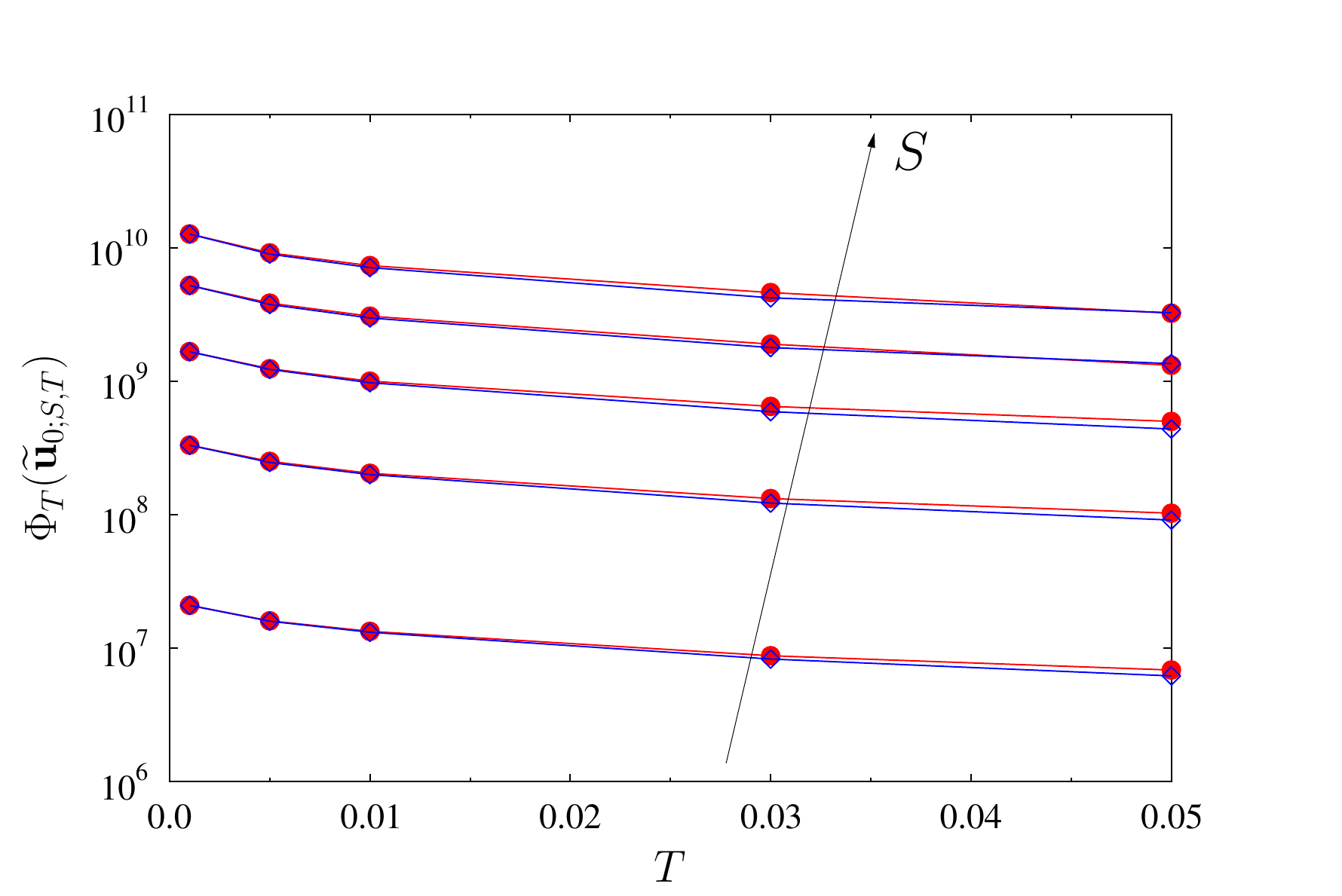}}}
\caption{Dependence of the maxima of the objective functionals (a)
  $\Phi_T(\tuBT)$ from Problem \ref{pb:PhiL4} with $B^4 = 1{,}000,
  3{,}000, 6{,}000, 9{,}000, 12{,}000$ and (b) $\Phi_T(\tuST)$ from
  Problem \ref{pb:PhiH3/4} with $S^2 = 400, 800, 1{,}200, 1{,}600,
  2{,}000$ on the length $T$ of the optimization window. In panel (a)
  the blue and red lines correspond to the partially symmetric and
  asymmetric branches, whereas in panel (b) they correspond to the
  symmetric and two-component branches. Arrows indicate the directions
  of increase of the constraints $B$ and $S$.}
\label{fig:maxLPSvsT}
\end{center}
\end{figure}

We now go on to discuss the structure of the extremal flows belonging
to the different maximizing branches by characterizing their symmetry
properties. We will do this by focusing on the componentwise
enstrophies \eqref{eq:Ei} whose time evolution in
representative solutions of Problems \ref{pb:PhiL4} and
\ref{pb:PhiH3/4} from both maximizing branches is shown in Figures
\ref{fig:EiL}a,b and \ref{fig:EiH}a,b, respectively.  As regards
solutions of Problem \ref{pb:PhiL4} corresponding to the dominating
branch which are shown in Figure \ref{fig:EiL}a, we have $\E_1(\u(t))
= \E_2(\u(t)) > \E_3(\u(t))$, $\forall t \in [0,T]$, indicating that
{in these flows} two vorticity components always contribute the same
amount of enstrophy.  On the other hand, for solutions corresponding
to the second branch, the {componentwise enstrophies} $\E_1(\u(t))$,
$\E_2(\u(t))$ and $\E_3(\u(t))$ remain distinct at almost all times $t
\in [0,T]$. We will thus refer to these two branches as ``partially
symmetric'' and ``asymmetric''. As concerns solutions of Problem
\ref{pb:PhiH3/4}, the results shown in Figures \ref{fig:EiH}a and
\ref{fig:EiH}b indicate that we have $\E_1(\u(t)) = \E_2(\u(t)) =
\E_3(\u(t))$ and $\E_1(\u(t)) = \E_2(\u(t)) > \E_3(\u(t)) = 0$,
$\forall t \in [0,T]$, for the two branches, which we will henceforth
refer to as ``symmetric'' and ``two-component'', respectively. In
solutions on these two branches the enstrophy is at all times
equipartitioned between two and three vorticity components.  We add
that these symmetry properties characterizing different branches are
robust and hold for different values of the parameters $B$, $S$ and
$T$.

\begin{figure}
\begin{center}
\mbox{\subfigure[]{\includegraphics[width=0.5\textwidth]{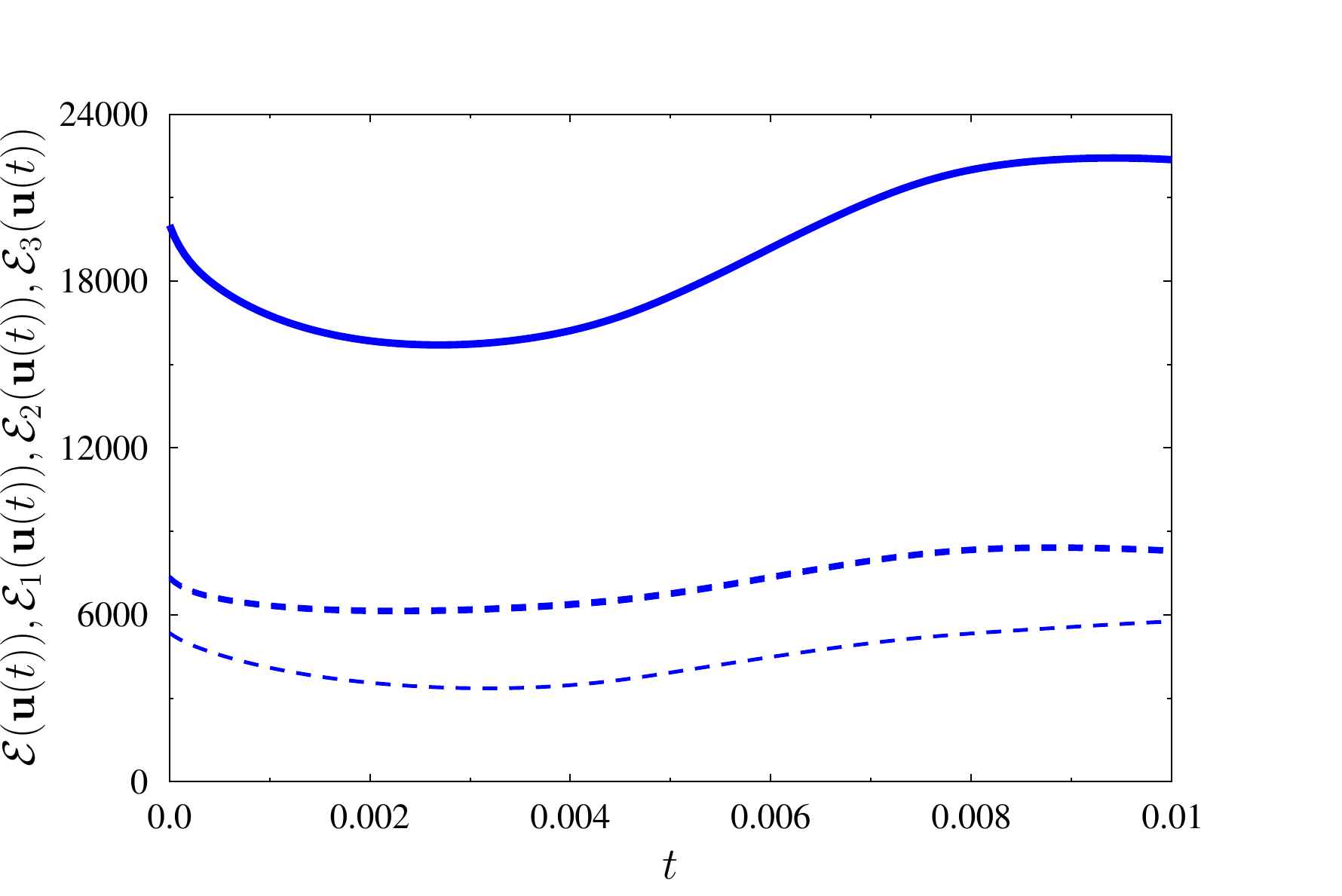}}\quad
\subfigure[]{\includegraphics[width=0.5\textwidth]{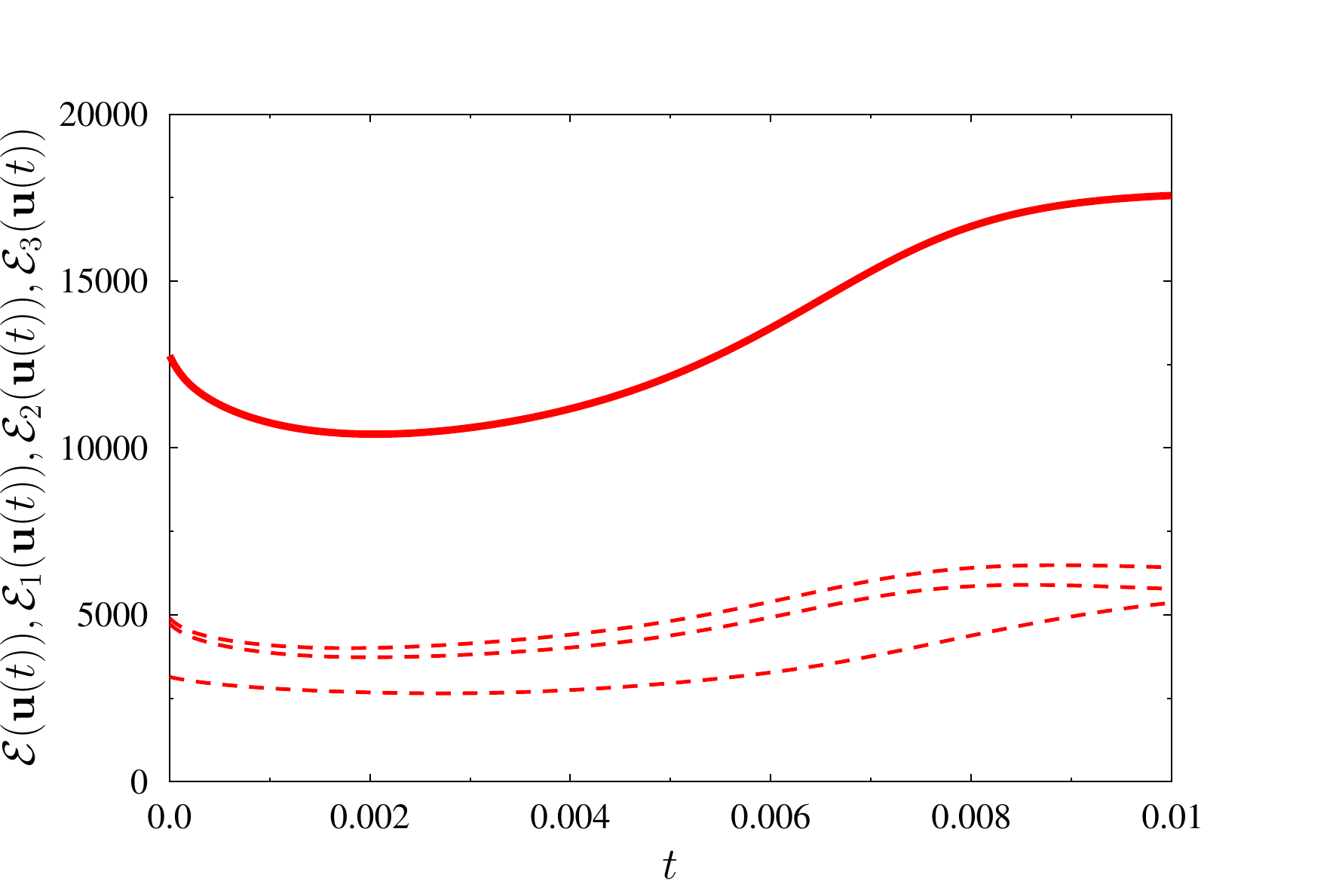}}}
\caption{Evolution of (thick solid lines) the total enstrophy
  $\E(\u(t))$ and (thin dashed lines) {the componentwise
    enstrophies} $\E_1(\u(t))$, $\E_2(\u(t))$, $\E_3(\u(t))$ in the
  solution of the Navier-Stokes system \eqref{eq:NS} with the optimal
  initial conditions $\tuBT$ on  (a) the partially symmetric
    branch and (b) the asymmetric branch obtained by solving Problem
    \ref{pb:PhiL4} with $B^4=12{,}000$ and $T=0.01$.}
\label{fig:EiL}
\mbox{\subfigure[]{\includegraphics[width=0.5\textwidth]{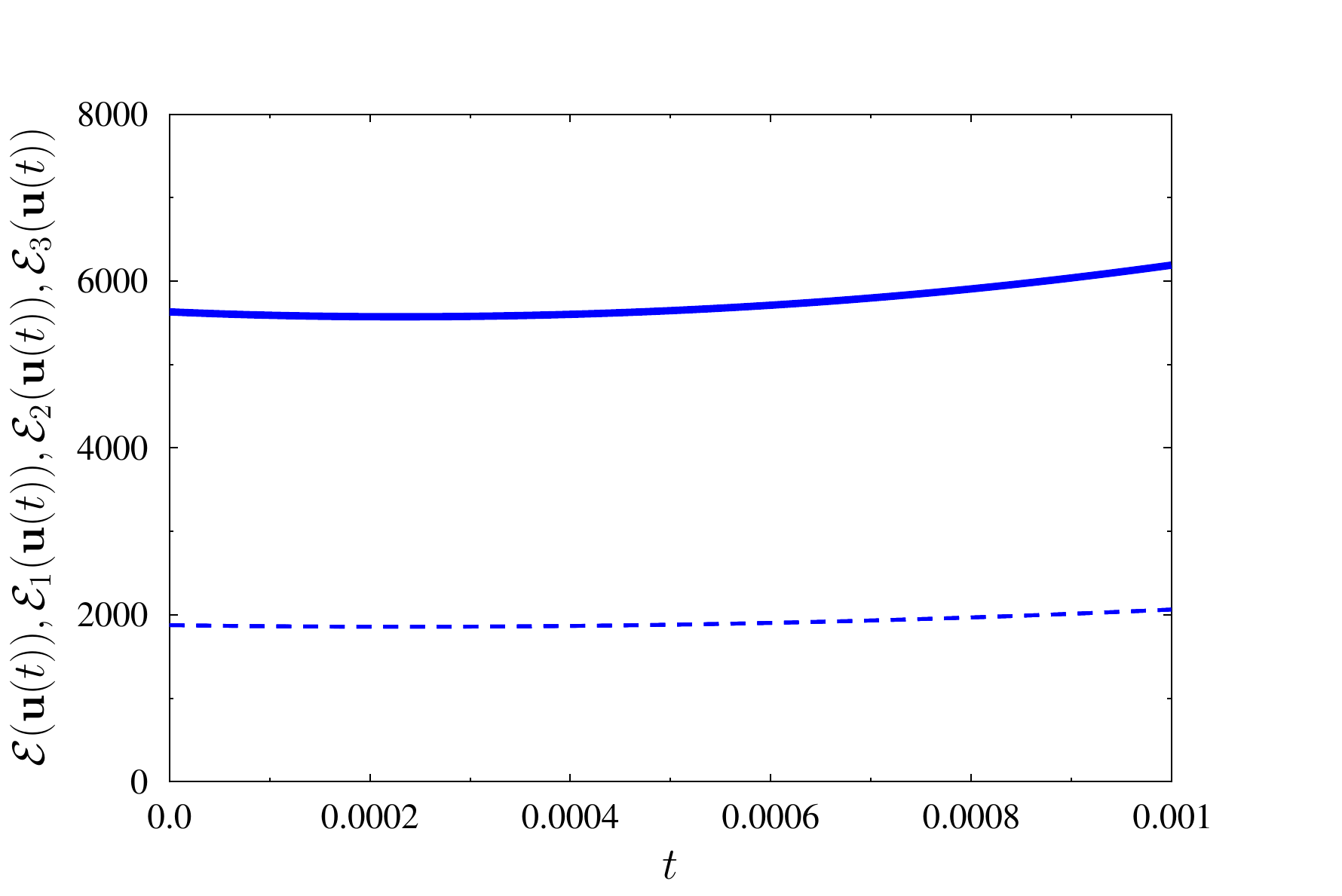}}\quad
\subfigure[]{\includegraphics[width=0.5\textwidth]{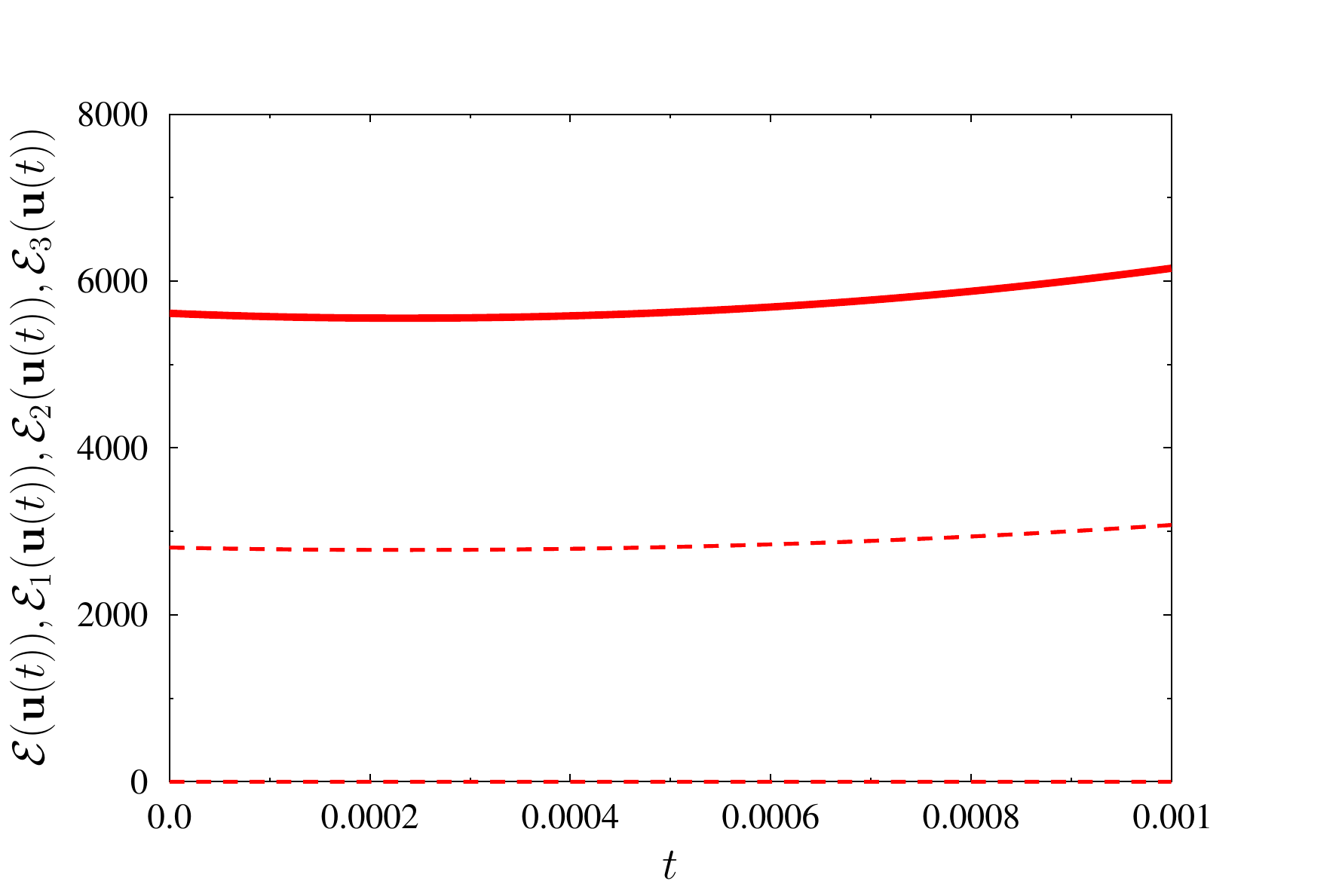}}}
\caption{Evolution of (thick solid lines) the total enstrophy
  $\E(\u(t))$ and (thin dashed lines) {the componentwise
    enstrophies} $\E_1(\u(t))$, $\E_2(\u(t))$, $\E_3(\u(t))$ in the
  solution of the Navier-Stokes system \eqref{eq:NS} with the optimal
  initial conditions $\tuST$ on (a) the symmetric branch and (b) the
  two-component branch obtained by solving Problem \ref{pb:PhiH3/4}
  with $S^2=2{,}000$ and $T=0.001$.}
\label{fig:EiH}
\end{center}
\end{figure}

In order to understand the physical structure of the extreme flows,
the optimal initial conditions $\tuBT$ and $\tuST$, obtained by
solving Problems \ref{pb:PhiL4} and \ref{pb:PhiH3/4} are shown in
Figures \ref{fig:tuBT}a,b and \ref{fig:tuST}a,b. In both cases they
were obtained with the largest considered values of the constraints,
i.e., $B^4 = 12{,}000$ and $S^2 = 2{,}000$. For Problem
\ref{pb:PhiL4}, the initial conditions shown were obtained with $T =
0.01$, which is the length of the time window for which the largest
value of the objective functional $\Phi_T(\tuBT)$ was attained,
cf.~Figure \ref{fig:maxLPSvsT}a. The initial condition corresponding
to the dominating partially-symmetric branch has the form of two
nearly parallel curved vortex {sheets}, cf.~Figure \ref{fig:tuBT}a. On
the other hand, the initial condition corresponding to the asymmetric
branch has the form of a single curved vortex sheet, cf.~Figure
\ref{fig:tuBT}b.  The time evolutions of the flows corresponding to
the optimal initial conditions shown in Figures \ref{fig:tuBT}a and
\ref{fig:tuBT}b are visualized in
\href{https://youtu.be/LHN_Z7orxEE}{Movie 1} and
\href{https://youtu.be/7jwciCRzbL4}{Movie 2} available on-line.  For
Problem \ref{pb:PhiH3/4} with the shortest considered time window $T =
0.001$ which also produced the largest value of the objective
functional $\Phi_T(\tuST)$, cf.~Figure \ref{fig:maxLPSvsT}b, in
Figures \ref{fig:tuST}a,b we see that the optimal initial condition
$\tuST$ is very similar for both branches and has the form of a single
vortex ring. The only difference is that the axis of the vortex ring
is aligned with one of the coordinate directions in the case of the
two-component branch and with the diagonal direction of the domain
$\Omega$ for the symmetric branch.  This property explains the
equipartition of enstrophy observed in Figures \ref{fig:EiH}a and
\ref{fig:EiH}b.  As the value of the constraint $S$ increases or the
time window $T$ shrinks, the corresponding optimal initial conditions
$\tuST$ become more localized such that the orientation of the vortex
structure with respect to the domain $\Omega$ plays a lesser role.
This explains why the optimal initial data from the two branches
obtained in Problem \ref{pb:PhiH3/4} yield very similar values of the
objective function $\Phi_T(\tuST)$, cf.~Figure \ref{fig:maxLPSvsT}b.
{The time evolution of the flow corresponding to the optimal initial
  condition shown in Figures \ref{fig:tuST}a is visualized in
  \href{https://youtu.be/4vrQqFZ5Hi8}{Movie 3} available on-line. We
  see that this evolution involves the translation and diffusion of
  the vortex ring.}

\begin{figure}
\begin{center}
\mbox{\subfigure[]{\includegraphics[width=0.45\textwidth]{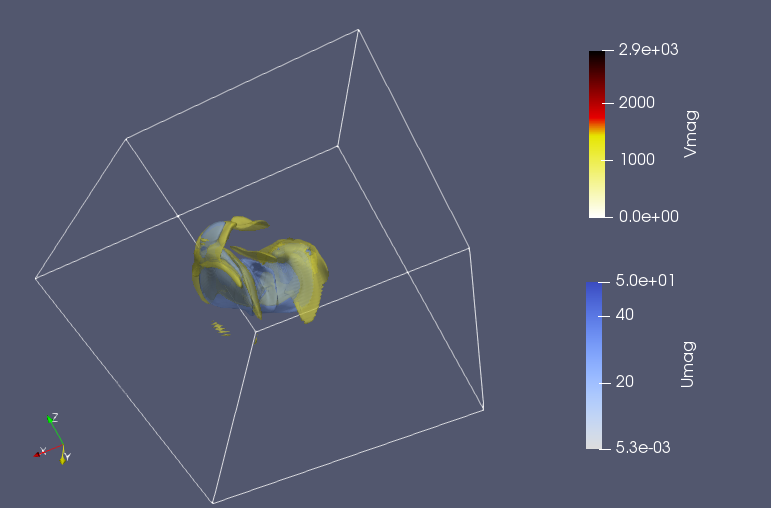}}\qquad
\subfigure[]{\includegraphics[width=0.45\textwidth]{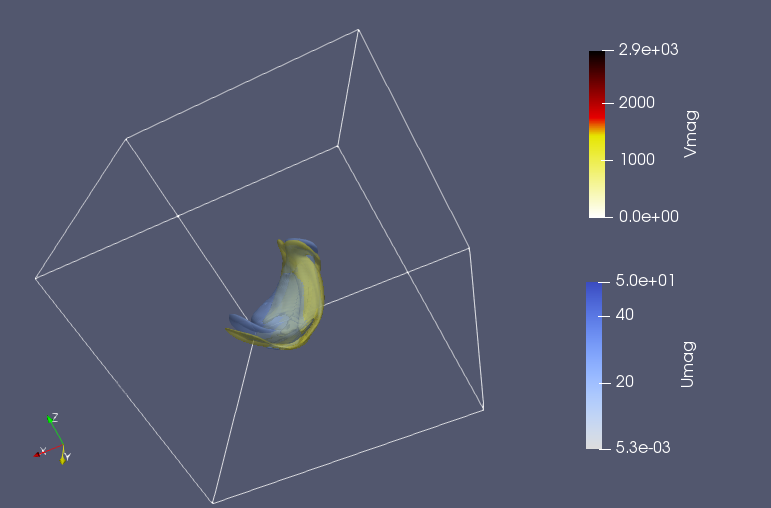}}}
\caption{Optimal initial conditions $\tuBT$ on (a) the partially
  symmetric branch and (b) asymmetric branch obtained by solving
  Problem \ref{pb:PhiL4} with $B^4=12{,}000$ and $T=0.01$.  Yellow and
  blue represent the iso-surfaces of the vorticity magnitude
  $\left|\left(\bnabla \times \tuBT\right)(\x)\right|$ and the
  velocity magnitude $\left|\tuBT(\x)\right|$, respectively. {The time
    evolutions of the flows corresponding these initial conditions are
    visualized in \href{https://youtu.be/LHN_Z7orxEE}{Movie 1} and
    \href{https://youtu.be/7jwciCRzbL4}{Movie 2} available as on-line.}}
\label{fig:tuBT}
\bigskip\bigskip
\mbox{\subfigure[]{\includegraphics[width=0.45\textwidth]{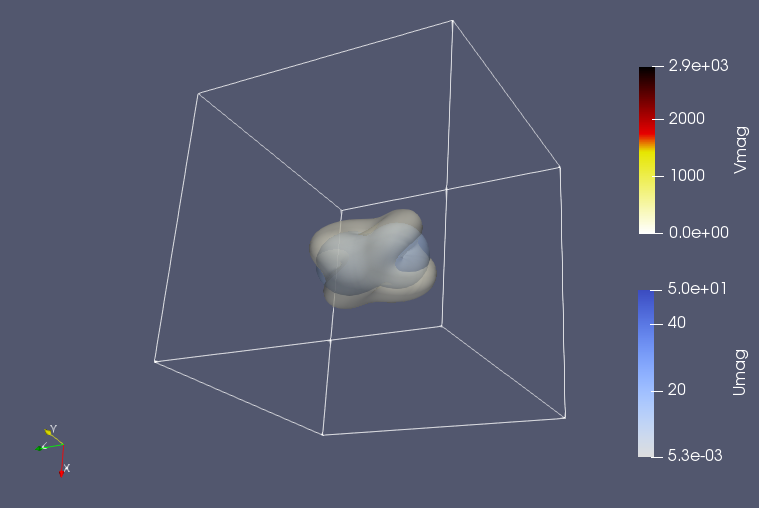}}\qquad
\subfigure[]{\includegraphics[width=0.45\textwidth]{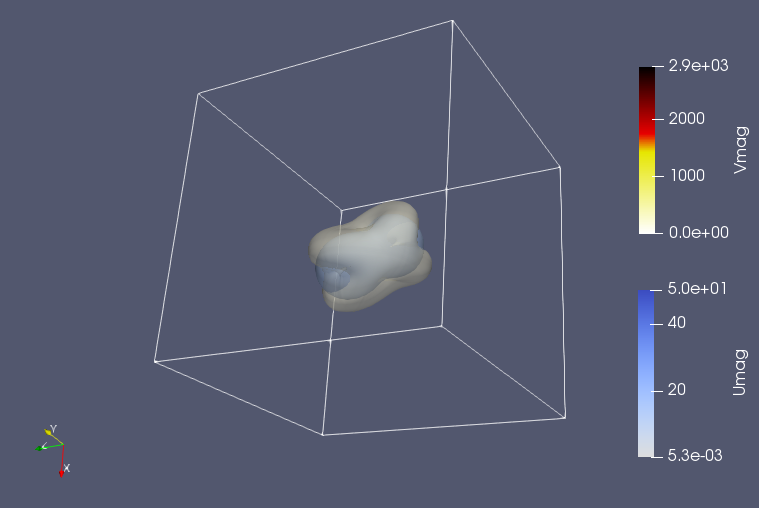}}}
\caption{Optimal initial conditions $\tuST$ on (a) the two-component
  branch and (b) symmetric branch obtained by solving Problem
  \ref{pb:PhiH3/4} with $S^2=2{,}000$ and $T=0.001$.  Yellow and blue
  represent the iso-surfaces of the vorticity magnitude
  $\left|\left(\bnabla \times \tuST\right)(\x)\right|$ and the
  velocity magnitude $\left|\tuST(\x)\right|$, respectively.  {The
    time evolution of the flow corresponding the initial condition
    shown in (a) is visualized in
    \href{https://youtu.be/4vrQqFZ5Hi8}{Movie 3} available on-line.}}
\label{fig:tuST}
\end{center}
\end{figure}

We now return to the question whether the quantity in
\eqref{eq:LPSblowup} with $q=4$ can become unbounded in finite time,
which would signal singularity formation. The results summarized in
Figures \ref{fig:maxLPSvsT}a and \ref{fig:maxLPSvsT}b show no evidence
of unbounded growth of the functional $\Phi_T(\u_0)$ when it is
maximized by solving Problems \ref{pb:PhiL4} and \ref{pb:PhiH3/4}. The
maximum growth achieved by this functional is presented in Figures
\ref{fig:maxLPS_vs_BS}a and \ref{fig:maxLPS_vs_BS}b where we plot
$\max_T \Phi_T(\tuBT)$ and $\max_T \Phi_T(\tuST)$, respectively, as
functions of the constraints $B$ and $S$. In other words, the maxima
are taken over a maximizing branch with a fixed value of the
constraint with respect to the length $T$ of the optimization window.
As is evident from Figures \ref{fig:maxLPS_vs_BS}a and
\ref{fig:maxLPS_vs_BS}b, both $\max_T \Phi_T(\tuBT)$ and $\max_T
\Phi_T(\tuST)$ reveal clear power-law dependence on the values of
the constraint which can be described by the following relations
obtained by performing least-squares fits
\begin{subequations}
\label{eq:fits}
\begin{align}
\max_T \Phi_T(\tuBT)  & \approx {(0.6478 \pm 0.1153)} \left(\| \tuBT \|_{L^4(\Omega)}^4 \right)^{2.261\pm 0.021}, \label{eq:fitB} \\
\max_T \Phi_T(\tuST)  & \approx {(9.308\pm 0.373)\times 10^{-4}} \left(\| \tuST \|_{\dot{H}^{3/4}(\Omega)}^2 \right)^{3.979\pm 0.005}. \label{eq:fitS}
\end{align}
\end{subequations}
From Figure \ref{fig:maxLPSvsT}b we conclude that in Problem 2, the
functional $\Phi_T(\tuST)$ achieves its maximum with respect to $T$ in
the limit $T \rightarrow 0$, and thus $\max_T\Phi_T(\tuST)$ depends
only on the value of the constraint $\| \tuST \|_{L^4(\Omega)}$.
Therefore, solving Problem 2 for $T \rightarrow 0$ is equivalent to
seeking a divergence-free vector field with a fixed $H^{\frac{3}{4}}$
seminorm and a maximum $L^4$ norm, which explains the presence of an
exponent close to 4 in \eqref{eq:fitS}.  As a result, the optimal
initial data $\tuST$ obtained for different values of the constraint
$S$ with $T \rightarrow 0$ are identical up to normalization.

\begin{figure}
\begin{center}
\mbox{\subfigure[]{\includegraphics[width=0.5\textwidth]{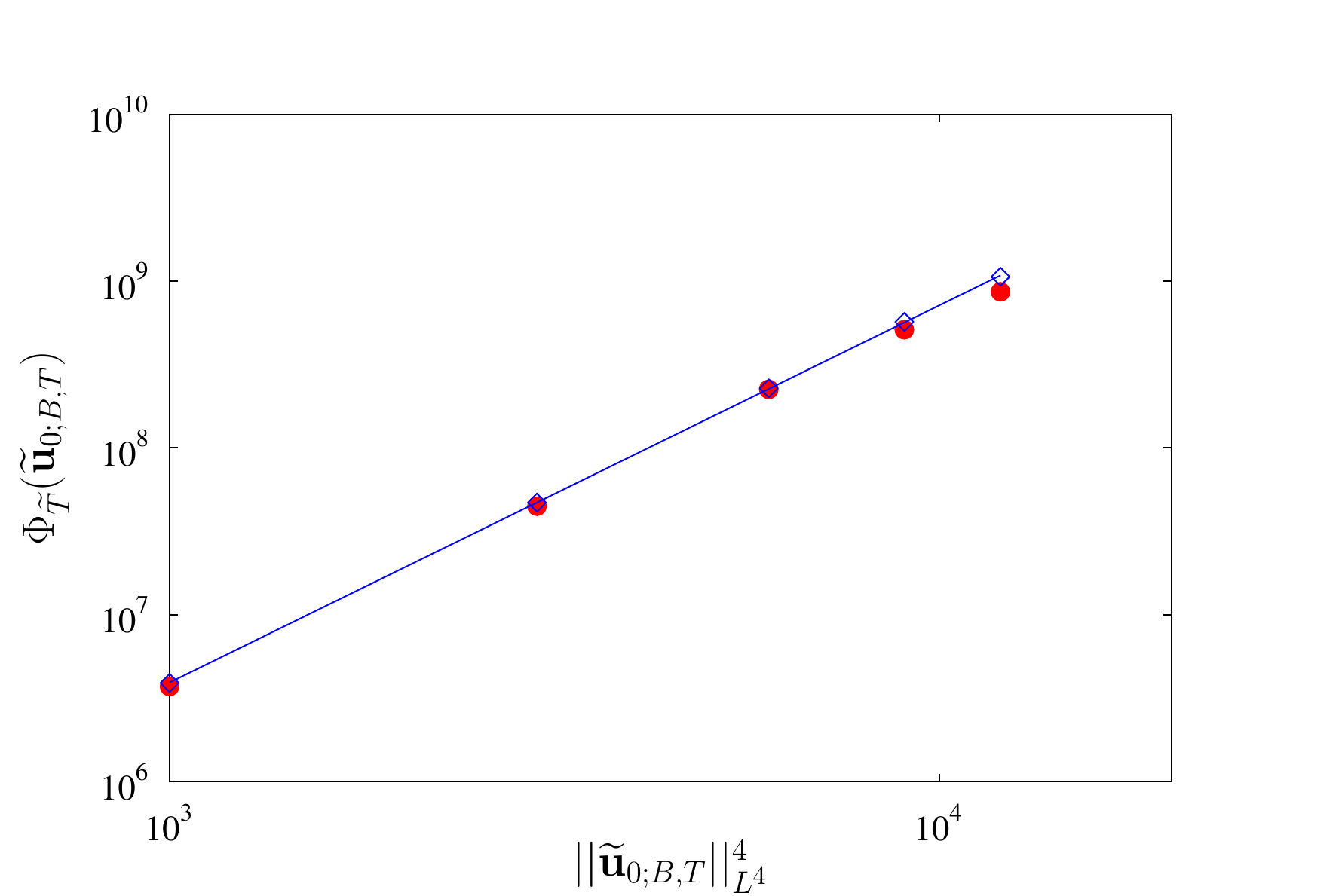}} \quad
\subfigure[]{\includegraphics[width=0.5\textwidth]{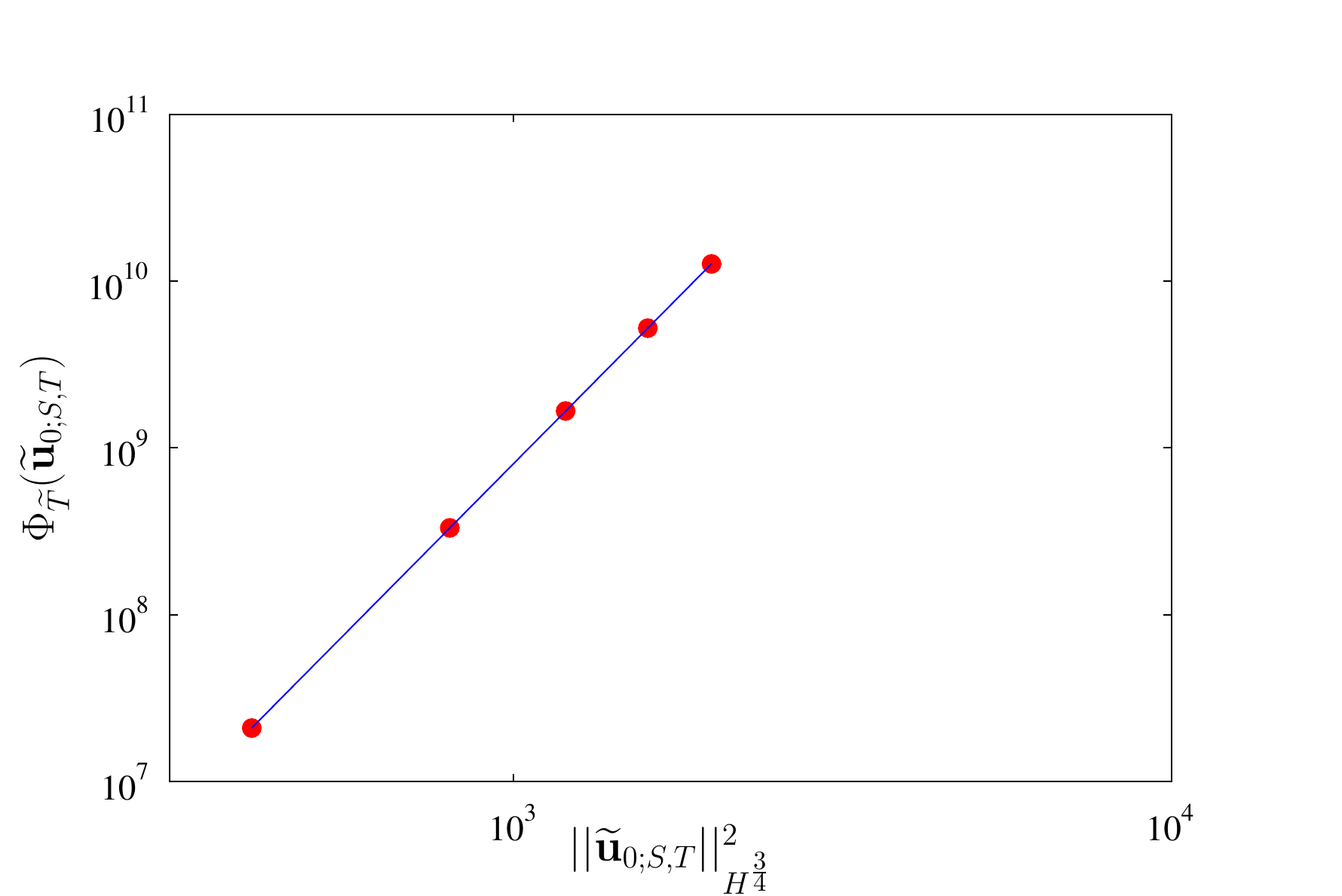}}}
\caption{Dependence of (a) $\max_T \Phi_T(\tuBT)$ on $B^4 = \|
  \tuBT\|_{L^4(\Omega)}^4$ and (b) $\max_T \Phi_T(\tuST)$ on $S^2 = \|
  \tuST\|_{\dot{H}^{3/4}(\Omega)}^2$ for Navier-Stokes flows with the
  optimal initial conditions $\tuBT$ and $\tuST$ obtained by solving
  Problems \ref{pb:PhiL4} and \ref{pb:PhiH3/4}, respectively. In panel
  (a) blue diamonds and red circles correspond to the partially
  symmetric and asymmetric branches, whereas in panel (b) these
  symbols correspond to the symmetric and two-component branches.
  Solid lines represent least-squares fits \eqref{eq:fitB} to the data
  from the partially symmetric branch in panel (a) and \eqref{eq:fitS}
  to the data from the symmetric branch in panel (b).}
\label{fig:maxLPS_vs_BS}
\end{center}
\end{figure}

The results obtained for Problem \ref{pb:PhiL4} can also provide
insights about the sharpness of the instantaneous estimate
\eqref{eq:dLqdt}. More specifically, as shown in Appendix
\ref{sec:dLqdt}, solutions of Problem \ref{pb:PhiL4} in the limit $T
\rightarrow 0$ approximate solutions of the instantaneous optimization
problem $\max_{\u \in {\L}_{B}} \frac{d}{dt}
\|\u\|_{L^4(\Omega)}^4$, where $\frac{d}{dt} \|\u\|_{L^4(\Omega)}$ can
be expressed using the Navier-Stokes system \eqref{eq:NS}.  Figure
\ref{fig:dLqdt} shows the dependence of $\frac{d}{dt}
\|\u(t)\|_{L^4(\Omega)}^4\big|_{t=0}$ approximated numerically based
on the solution of the Navier-Stokes system \eqref{eq:NS} with the
optimal initial condition $\tuBT$ obtained from Problem \ref{pb:PhiL4}
with the shortest considered optimization window $T = 0.001$ on
$\|\tuBT\|_{L^4(\Omega)}^4$.  For both branches the figure reveals an
essentially the same power-law relation
\begin{equation}
\frac{d}{dt} \|\u(t)\|_{L^4(\Omega)}^4\big|_{t=0} \approx (221.5\pm 104) \left( \|\tuBT)\|_{L^4(\Omega)}^4 \right)^{1.117\pm 0.05}.
\label{eq:fit_dLqdt}
\end{equation}
It is clear that the exponent 1.117 in \eqref{eq:fit_dLqdt} is
significantly smaller than the exponent of 3 predicted by estimate
\eqref{eq:dLqdt} with $q = 4$.

\begin{figure}
\begin{center}
\mbox{\subfigure[]{\includegraphics[width=0.5\textwidth]{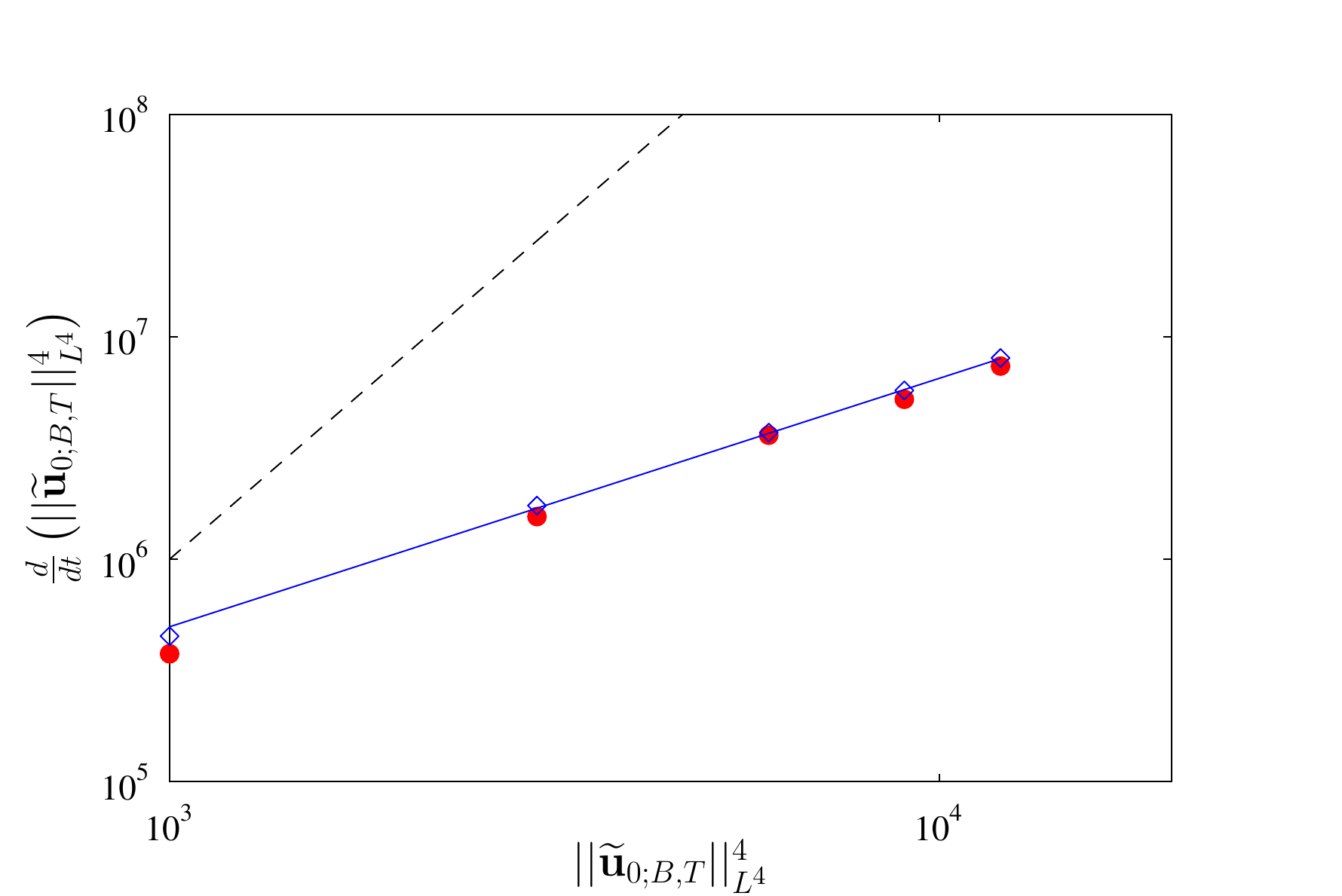}}}
\caption{Dependence of the rate of change $\frac{d}{dt}
  \|\u(t)\|_{L^4(\Omega)}^4\big|_{t=0}$ on $\|\tuBT\|_{L^4(\Omega)}^4$
  for Navier-Stokes flows with the optimal initial conditions $\tuBT$
  obtained by solving Problems \ref{pb:PhiL4} with $T = 0.001$, which
  is the shortest considered time window. The blue diamonds and red
  circles correspond to the partially symmetric and asymmetric
  branches, respectively, whereas the blue solid line represents the
  least-squares fit \eqref{eq:fit_dLqdt}.  The dashed black line
  corresponds to the exponent of 3 obtained in \eqref{eq:dLqdt} with
  $q=4$.}
\label{fig:dLqdt}
\end{center}
\end{figure}

Finally, we compare the extreme flows analyzed above to the extreme
flows constructed in \cite{KangYumProtas2020} in terms of the relative
growth of enstrophy.  Dependence of the maximum attained enstrophy
$\max_{t \ge 0} \E(t)$ on the initial enstrophy $\E_0$ in
Navier-Stokes flows with the optimal initial conditions $\tuBT$ and
$\tuST$ obtained by solving Problems \ref{pb:PhiL4} and
\ref{pb:PhiH3/4} is shown in Figure \ref{fig:maxE}, where, for
comparison, we also show relation \eqref{eq:maxT_vs_E0} discovered in
\cite{KangYumProtas2020}. The corresponding least-squares fits have
the form
\begin{subequations}
\begin{align}
\max_{t \ge 0} \E (t) & \approx  (1.697\pm 3.776)\times 10^{-2} \, \E_0 ^{1.499\pm 0.249} , \label{eq:fit_maxET_L} \\
\max_{t \ge 0} \E (t) & \approx (0.5155\pm 0.0532) \, \E_0 ^{1.147\pm 0.013}. \label{eq:fit_maxET_H}
\end{align}
\label{eq:fit_maxET}
\end{subequations}
It is intriguing to note that that the power-law relation
\eqref{eq:fit_maxET_L} corresponding to the partially-symmetric branch
obtained in Problem \ref{pb:PhiL4} features an essentially the same
exponent close to 3/2 as in \eqref{eq:maxT_vs_E0}, although the
prefactor is much smaller.

\begin{figure}
\begin{center}
\mbox{\subfigure[]{\includegraphics[width=0.5\textwidth]{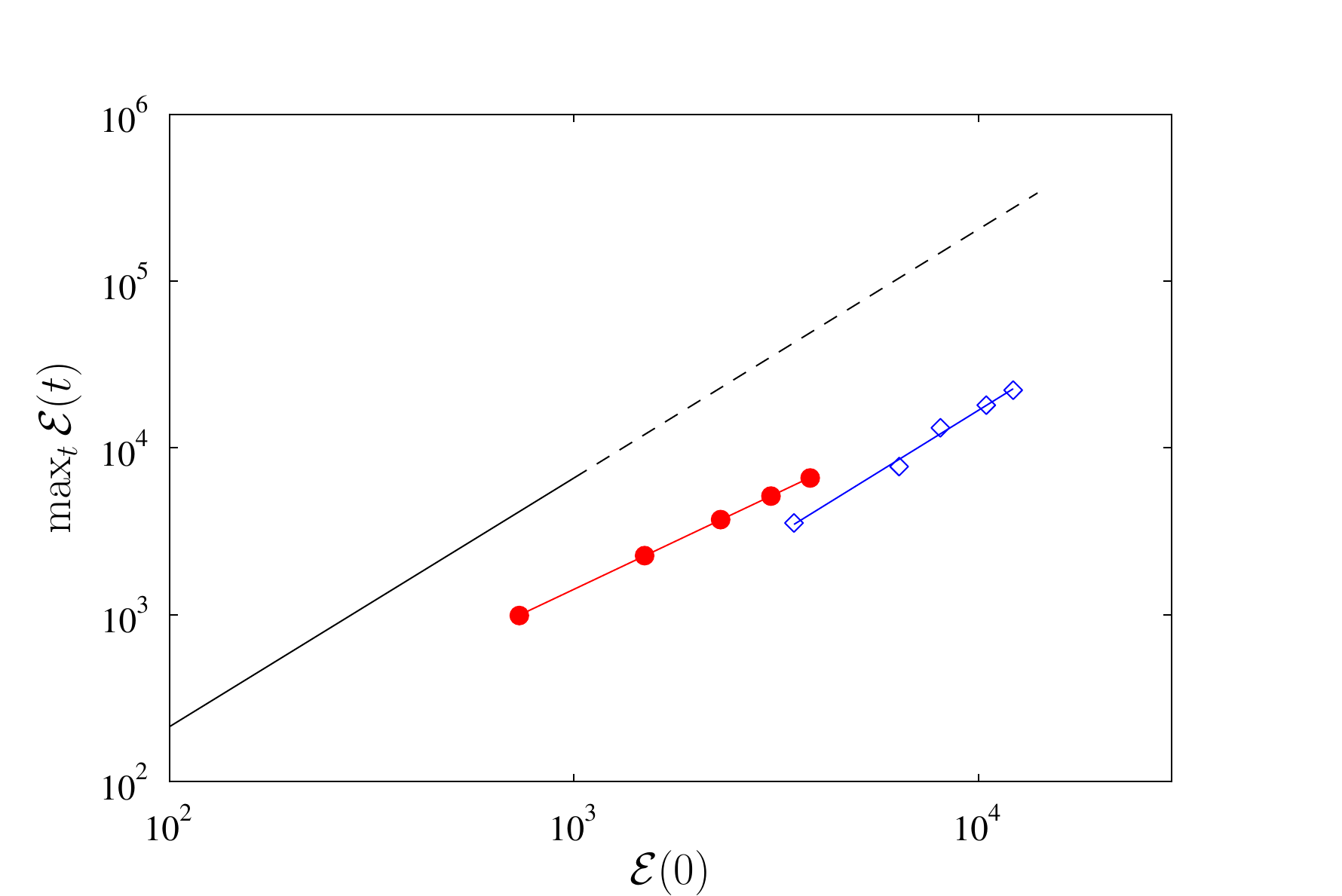}}} 
\caption{Dependence of the maximum attained enstrophy $\max_{t \ge 0}
  \E(t)$ on the initial enstrophy $\E_0$ in Navier-Stokes flows with
  the optimal initial conditions (blue diamonds) $\tuBT$ obtained by
  solving Problems \ref{pb:PhiL4} (the partially-symmetric branch) and
  (red circles) $\tuST$ obtained by solving Problems \ref{pb:PhiH3/4}
  (the symmetric branch). Each symbol corresponds to a different value
  of the constraint $B$ and $S$, and in all cases the results are
  presented for the length $T$ of the optimization window producing
  the largest value of $\max_{t \ge 0} \E(t)$. The solid black line
  represents the relation $\max_{t \ge 0} \E(t) = C \E_0^{3/2}$
  discovered in \cite{KangYumProtas2020} and the dashed solid line its
  extrapolation to higher enstrophy values.}
\label{fig:maxE}
\end{center}
\end{figure}

\FloatBarrier

\subsection{Flows Obtained as Solutions of Problem \ref{pb:PsiL2}}
\label{sec:resultsP3}

Solution of Problem \ref{pb:PsiL2} has yielded a single maximizing
branch for each value of $\K_0$ with representative solutions shown in
Figure \ref{fig:P3L4} in terms of the time evolution of the norm
$\|\u(t)\|_{L^4(\Omega)}^4$ for ``short'' and ``long'' optimization
windows $T$. We see that, similarly to the solution of Problem
\ref{pb:PhiH3/4} in Figure \ref{fig:L4vst}b, the norm
$\|\u(t)\|_{L^4(\Omega)}^4$ is a decreasing function of time $t$. The
maximizing branches obtained for different values of the constraint
$\K_0$ are presented in terms of the dependence of the quantity
$T\Psi_T(\tuKT)$, which appears on the LHS of estimate
\eqref{eq:LPSbound}, on $T$ in Figure \ref{fig:P3vsT}. We see that for
each value of $\K_0$ the quantity $T\Psi_T(\tuKT)$ is an increasing
function of the length $T$ of the optimization window approaching a
certain limit as $T \rightarrow \infty$. In order to quantify its
behavior in this limit, for each value of $\K_0$ we construct a fit to
the dependence of $T\Psi_T(\tuKT)$ on $T$ in the form
\begin{equation}
g(T) := \psi_{\K_0} - \alpha \, e^{-\beta T}, \qquad T > 0,
\label{eq:g}
\end{equation}
where $\psi_{\K_0},\alpha,\beta \in \RR^+$ are parameters determined
via least-squares minimization, such that $\psi_{\K_0} \approx \lim_{T
  \rightarrow \infty} T\Psi_T(\tuKT)$.

\begin{figure}
\begin{center}
\mbox{\subfigure[]{\includegraphics[width=0.5\textwidth]{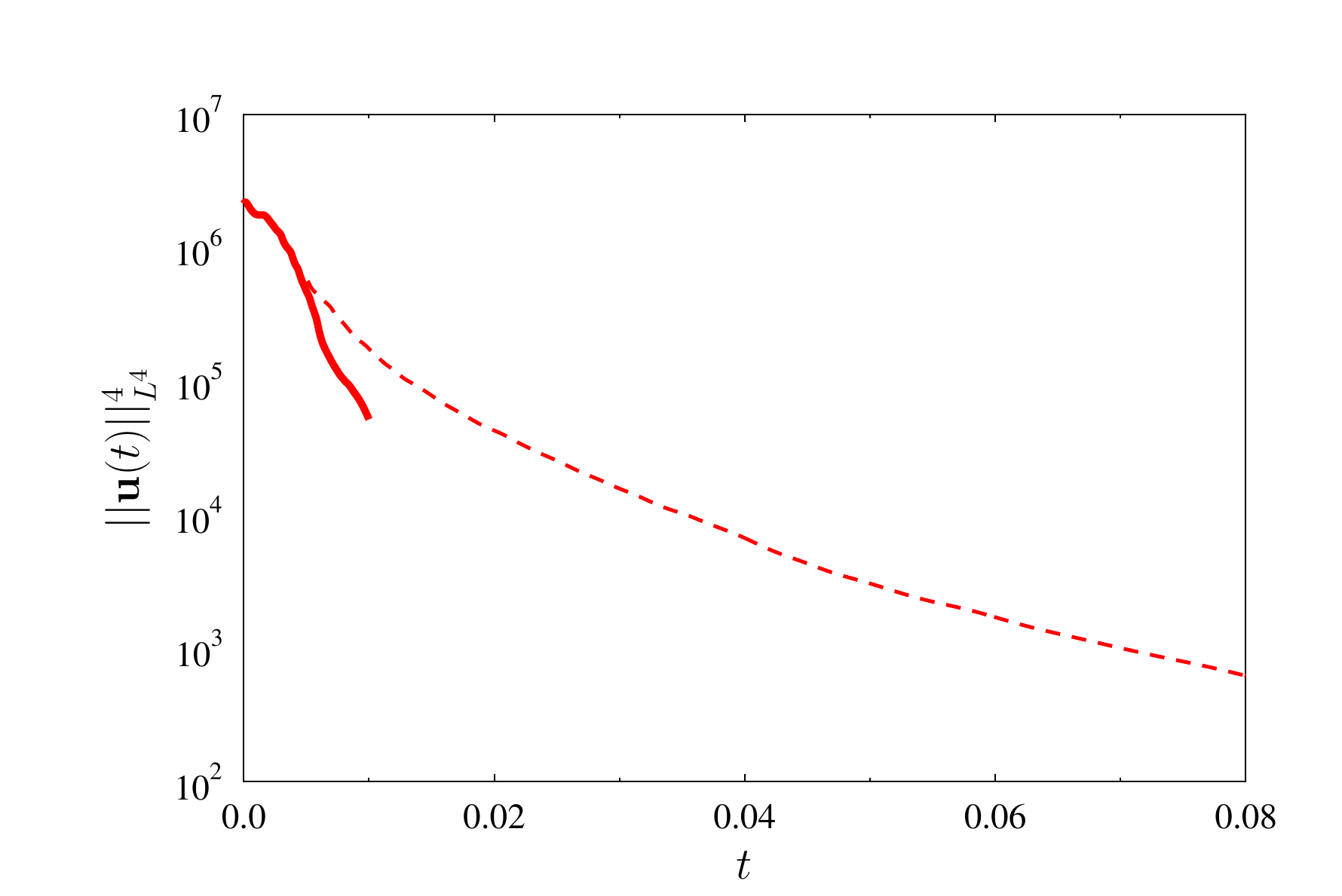}}} 
\caption{Time evolution of $\| \u (t) \|_{L^4}^{4}$ in the
  Navier-Stokes flows with optimal initial conditions $\tuKT$ obtained
  by solving Problem \ref{pb:PsiL2} with $\K_0=40$ using two different
  lengths $T=0.01$ and $T=0.08$ of the optimization window (in each
  case the results are shown on the time interval $[0,T]$ where
  optimization was performed).}
\label{fig:P3L4}
\end{center}
\end{figure}

\begin{figure}
\begin{center}
\mbox{\subfigure[]{\includegraphics[width=0.5\textwidth]{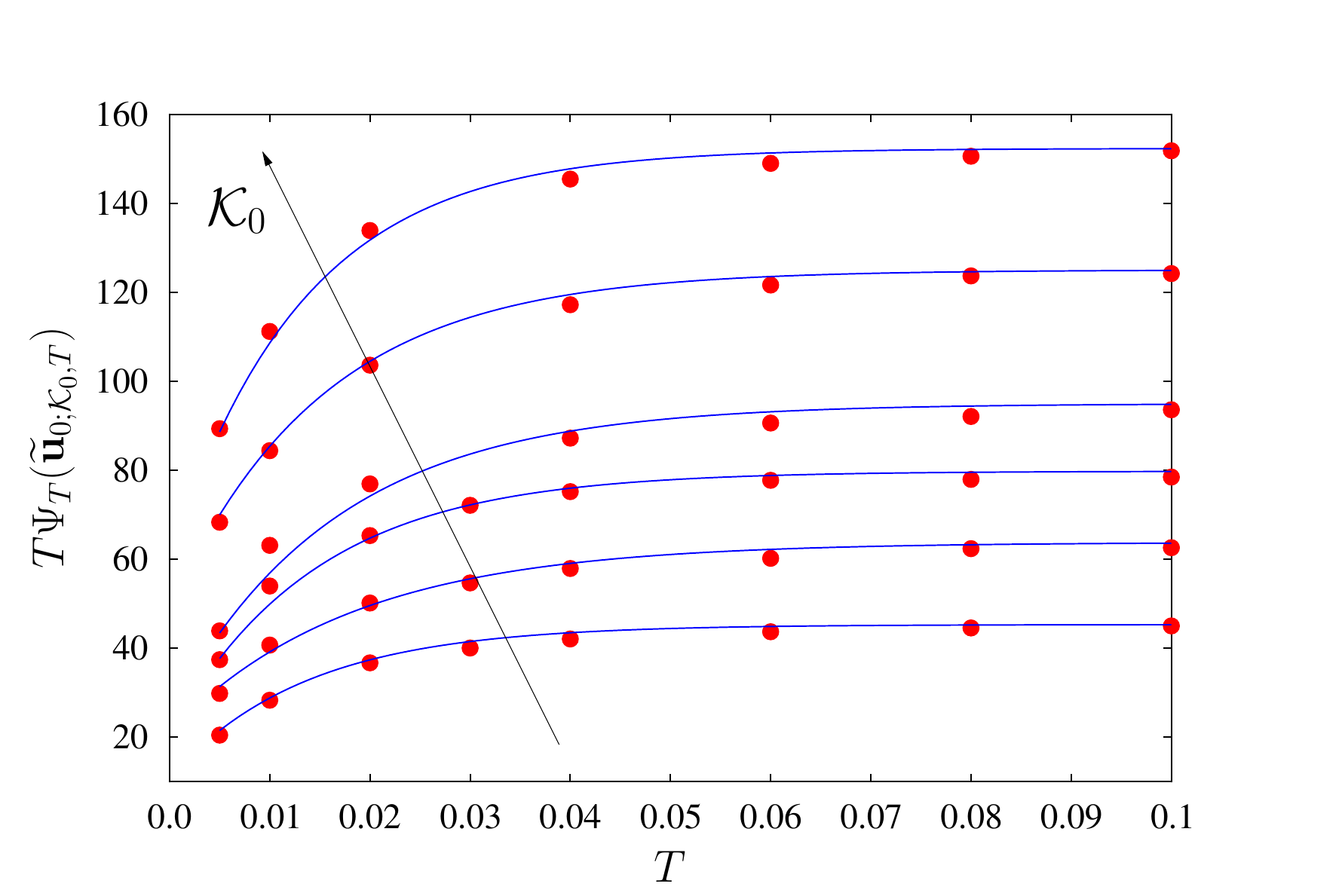}}} 
\caption{Dependence of the maxima of the quantity $T \Psi_T(\tuKT)$ on
  the length $T$ of the optimization window. Red symbols represent
  solutions of Problem \ref{pb:PsiL2} for $\K_0 = 15, 20, 25, 30, 40,
  50$ (arrow represents the direction of increase of $\K_0$) whereas
  blue lines are the fits obtained with formula \eqref{eq:g}.}
\label{fig:P3vsT}
\end{center}
\end{figure}

We now go on to discuss the structure of the extremal flows on the
maximizing branches by characterizing their symmetry properties using
the componentwise enstrophies \eqref{eq:Ei}. Their time evolution in
representative solutions of Problem \ref{pb:PsiL2} is shown in Figures
\ref{fig:EnsP3}a,b for short and long optimization windows $T$.  We
note that for both time windows we have the property $\E_1(\u(t)) =
\E_2(\u(t)) > \E_3(\u(t))$, $\forall t \in [0,T]$, the same as was
observed for solutions of Problem \ref{pb:PhiL4} on the dominating
branch, cf.~Figure \ref{fig:EiL}a. Hence, these optimal solutions can
be described as partially symmetric. However, in contrast to solutions
of Problem \ref{pb:PhiL4}, the time evolution of the enstrophy in
solutions of Problem \ref{pb:PsiL2} is much less regular and involves
significantly higher values. This more ``turbulent'' nature of
solutions of Problem \ref{pb:PsiL2} is also evident in the form of the
corresponding optimal initial conditions $\tuKT$ shown for the two
time windows in Figures \ref{fig:tuKT}a,b. As we can see, these
optimal initial conditions are less regular and involve more
small-scale features than the optimal initial condition obtained by
solving Problem \ref{pb:PhiL4}, cf.~Figure \ref{fig:tuBT}a. {The time
  evolution of the flow corresponding to the optimal initial condition
  shown in Figures \ref{fig:tuKT}a is visualized in
  \href{https://youtu.be/d0svCYCeeQQ}{Movie 4} available on-line. We
  see that this evolution involves the translation of a turbulent spot
  followed by its eventual bursting.}

\begin{figure}
\begin{center}
\mbox{\subfigure[]{\includegraphics[width=0.5\textwidth]{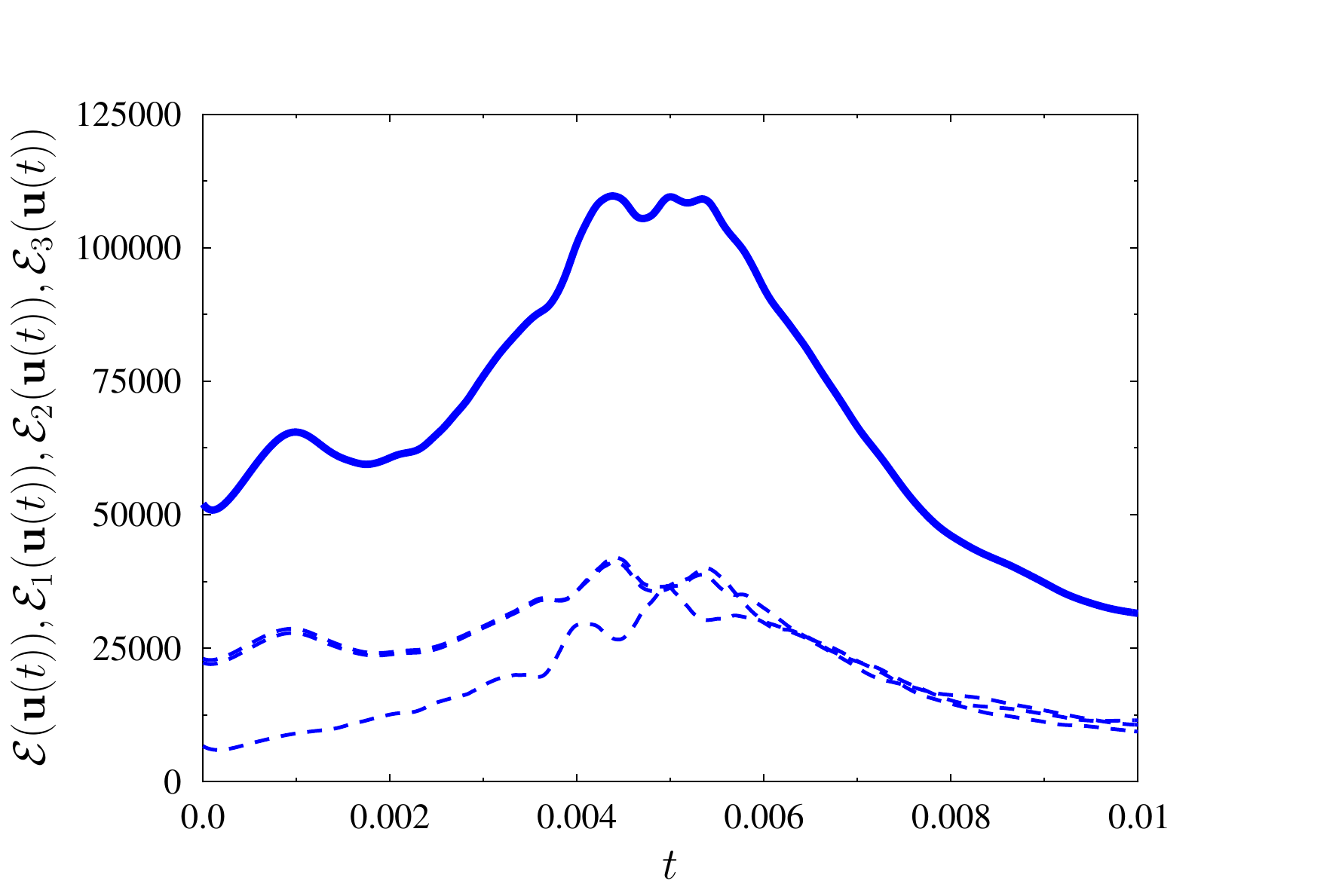}}\quad
\subfigure[]{\includegraphics[width=0.5\textwidth]{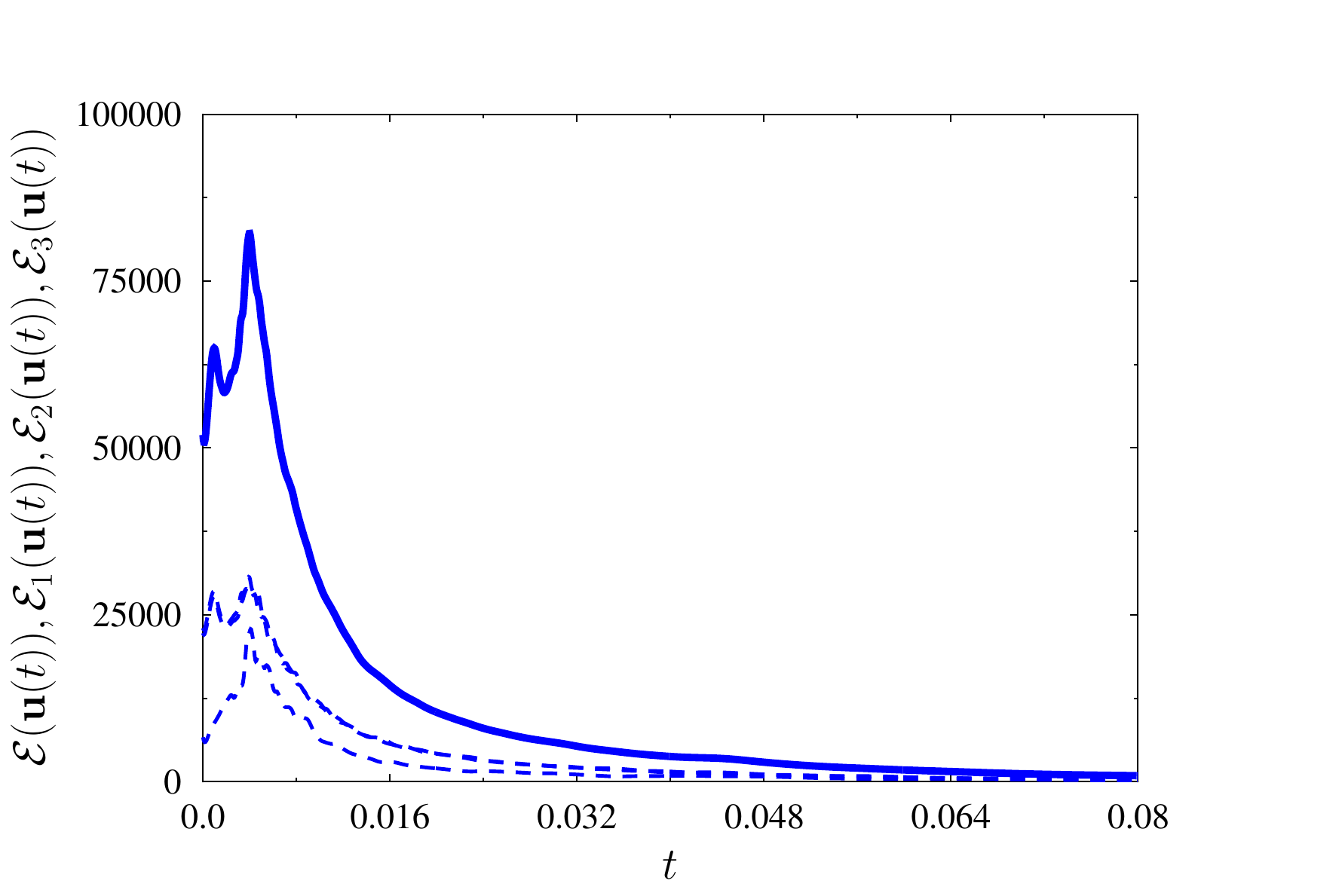}}}
\caption{Evolution of (thick solid lines) the total enstrophy
  $\E(\u(t))$ and (thin dashed lines) {the componentwise enstrophies}
  $\E_1(\u(t))$, $\E_2(\u(t))$, $\E_3(\u(t))$ in the solution of the
  Navier-Stokes system \eqref{eq:NS} with the optimal initial
  conditions $\tuKT$ obtained by solving Problem \ref{pb:PsiL2} with
  $\K_0=40$ and (a) $T=0.01$ and (b) $0.08$.}
\label{fig:EnsP3}
\bigskip\bigskip
\mbox{\subfigure[]{\includegraphics[width=0.45\textwidth]{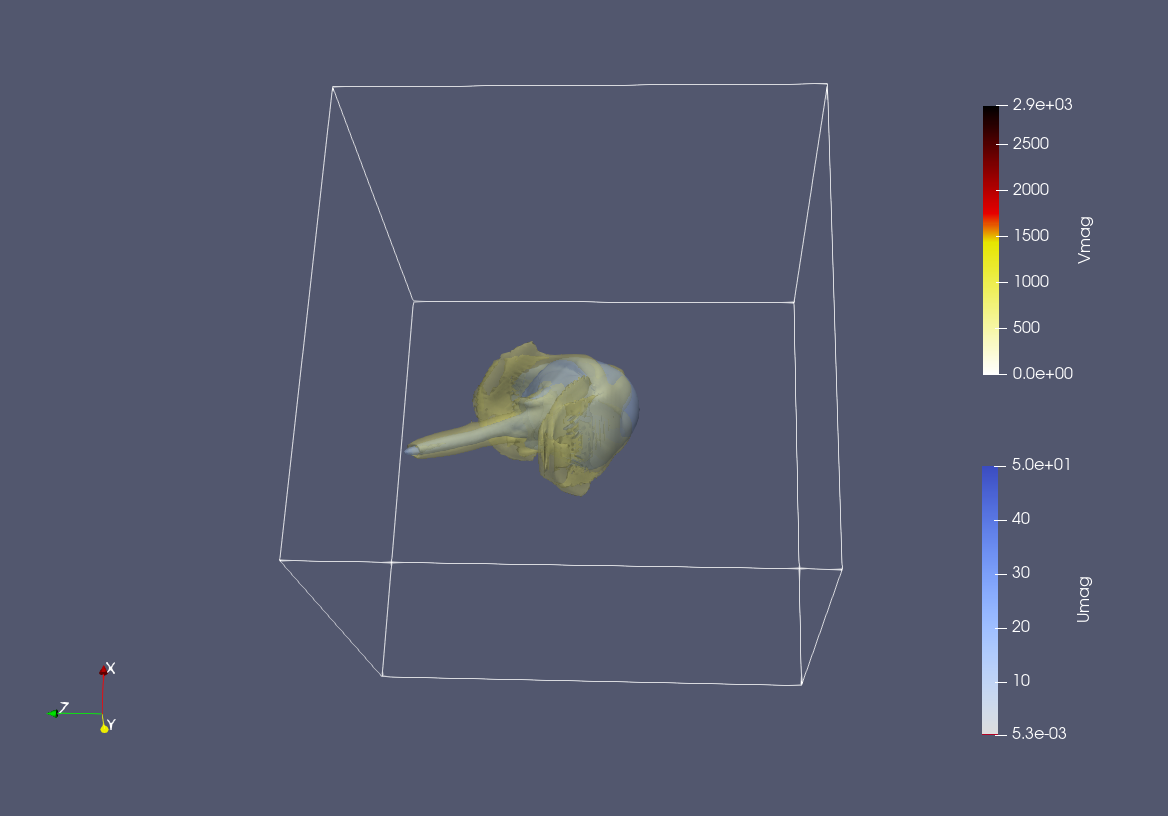}}\qquad
\subfigure[]{\includegraphics[width=0.45\textwidth]{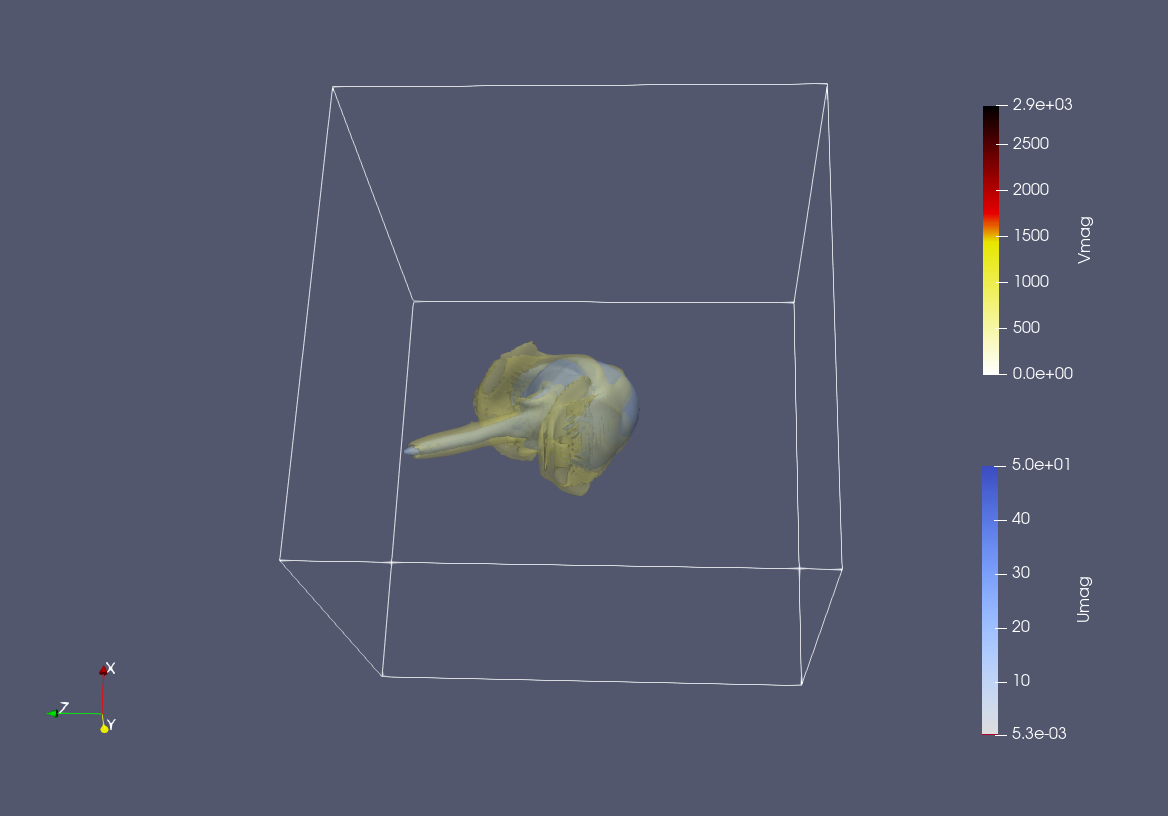}}}
\caption{Optimal initial conditions $\tuKT$ obtained by solving
  Problem \ref{pb:PsiL2} with $\K_0=40$ and (a) $T=0.01$ and (b)
  $T=0.08$.  Yellow and blue represent the iso-surfaces of the
  vorticity magnitude $\left|\left(\bnabla \times
      \tuKT\right)(\x)\right|$ and the velocity magnitude
  $\left|\tuKT(\x)\right|$, respectively. {The time evolution of the
    flow corresponding the initial condition shown in (a) is
    visualized in \href{https://youtu.be/d0svCYCeeQQ}{Movie 4}
    available on-line.}}
\label{fig:tuKT}
\end{center}
\end{figure}

Next, we analyze estimate \eqref{eq:LPSbound} in the limit of
long optimization windows $T$ where the term
$\|\u(T)\|_{L^2(\Omega)}^2$, cf.~\eqref{eq:LPSbound2a}, becomes
insignificant. To this end in Figure \ref{fig:P3vsK}a we plot
$\psi_{\K_0}$ from \eqref{eq:g} as function of $\K_0$ and
observe that
\begin{equation}
\psi_{\K_0} = \lim_{T   \rightarrow \infty} T\Psi_T(\tuKT) \approx (3.17\pm 0.9) \K_0^{0.998 \pm 0.082},
\label{eq:PsiK0}
\end{equation}
which reveals a power-law dependence on $\K_0$ although the range of
this quantity is not very extensive. The exponent is close to 1 which
is smaller than the exponent 4/3 predicted by estimate
\eqref{eq:LPSbound} with $q = 4$.

\begin{figure}
\begin{center}
\mbox{\subfigure[]{\includegraphics[width=0.5\textwidth]{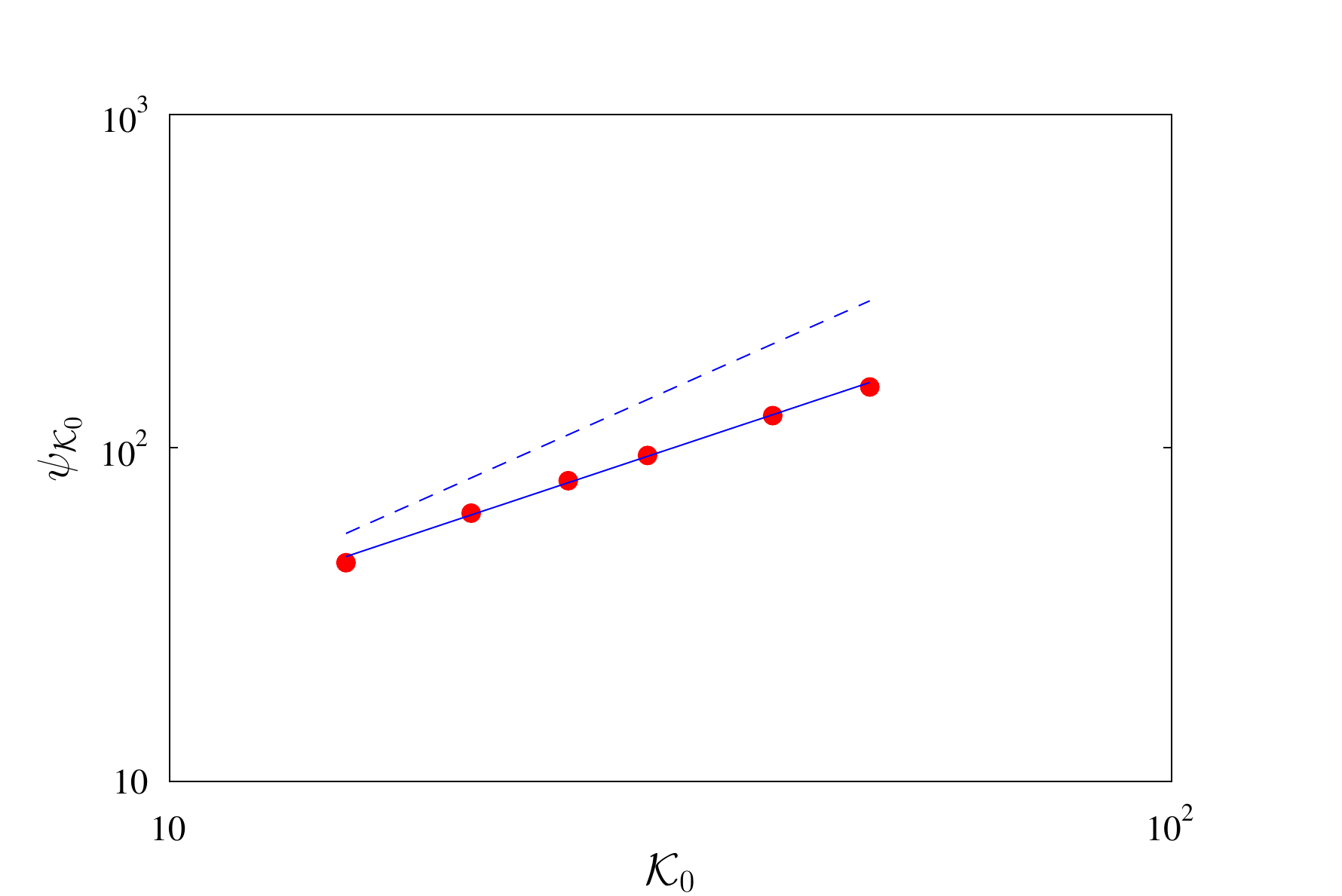}}\quad
\subfigure[]{\includegraphics[width=0.5\textwidth]{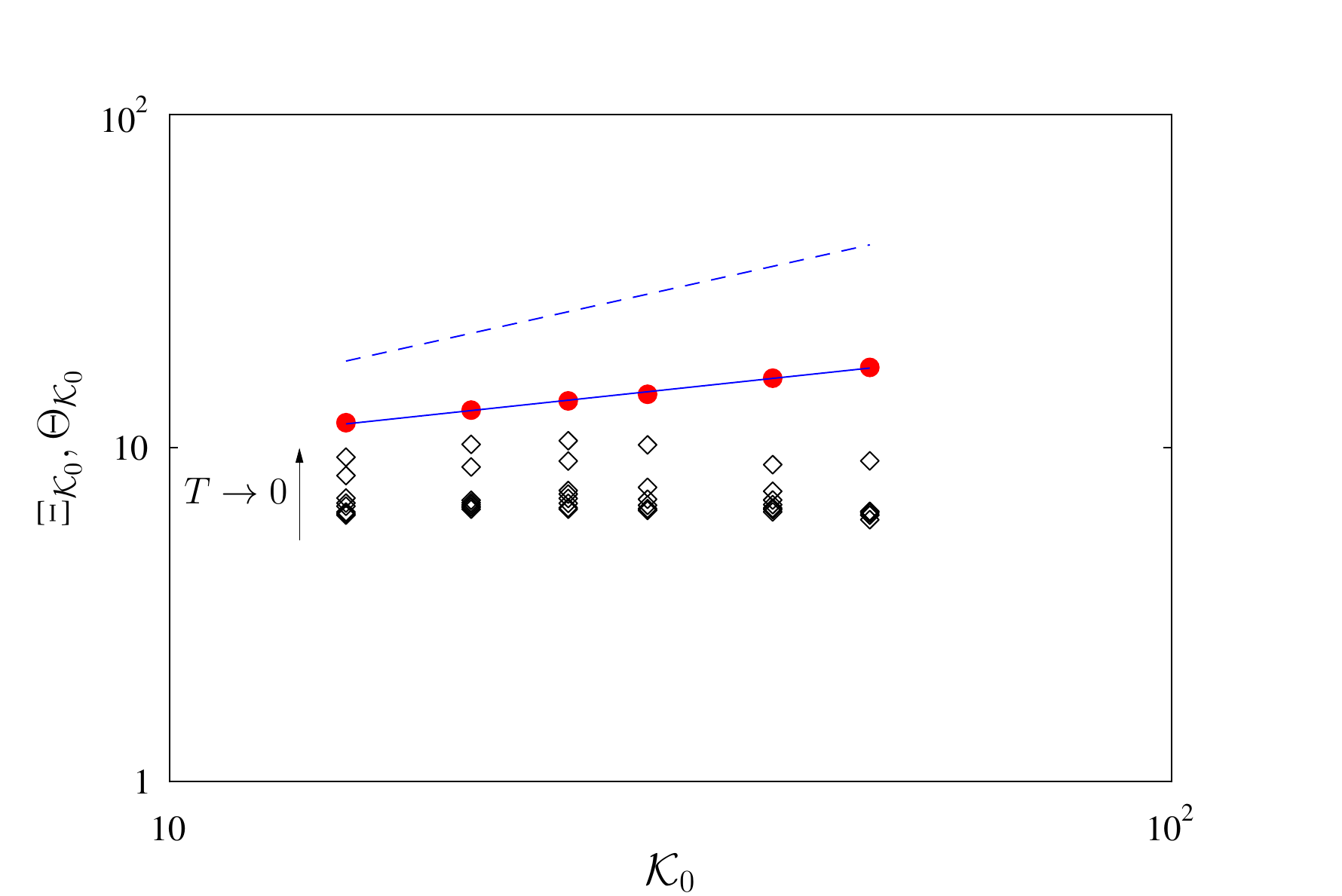}}} 
\caption{Dependence of (a) $\psi_{\K_0} \approx \lim_{T\rightarrow
    \infty} \Psi_T(\tuKT)$ and (b) (black diamonds) $\Xi_{\K_0}(T)$
  and (red circles) $\Theta_{\K_0}$ on $\K_0$ for Navier-Stokes flows
  with the optimal initial conditions $\tuKT$ obtained by solving
  Problems \ref{pb:PsiL2}. In (a) the dashed and solid lines
  represent, respectively, the expression on the RHS in estimate
  \eqref{eq:LPSbound} (with $q = 4$ and an arbitrarily chosen constant
  $C$) and the least-squares fit \eqref{eq:PsiK0}.  In (b) the dashed
  and solid lines represent, respectively, the expression on the RHS
  in \eqref{eq:Xi} (with an arbitrarily chosen constant $C$) and the
  least-squares fit \eqref{eq:ThetaK0}. The arrow indicates the trend
  with the {decrease} of $T$.}
\label{fig:P3vsK}
\end{center}
\end{figure}

In order to obtain insights about the properties of estimate
\eqref{eq:LPSbound} for short and intermediate times $T$, we consider
the relation
\begin{equation}
\Xi_{\K_0}(T) := \frac{\int_{0}^{T}||\u(t)||_{L^4}^{8/3} \, dt}{2\K_0  - \|\u(T)\|_{L^2(\Omega)}^2} \le \frac{C}{2^{1/3} \nu} \K_0^{2/3}
\label{eq:Xi}
\end{equation}
obtained by dividing \eqref{eq:LPSbound2} by the expression in
parentheses in \eqref{eq:LPSbound2b}, such that the dependence on $T$
is confined to the LHS. The quantity $\Xi_{\K_0}(T)$ is plotted as a
function of $\K_0$ for different time windows $T$ in Figure
\ref{fig:P3vsK}b. As we see in this figure, for a fixed $\K_0$,
$\Xi_{\K_0}(T)$ increases as $T$ is reduced. In order to better
understand the behavior of $\Xi_{\K_0}(T)$ for short optimization
windows we consider the limit $T \rightarrow 0$ and define
\begin{equation}
\Theta_{\K_0} := \lim_{T \rightarrow 0} \Xi_{\K_0}(T) 
=  \lim_{T \rightarrow 0} \frac{\Psi_T(\tuKT)}{2 \frac{\K_0 - \K(\u(T))}{T}} 
= \frac{\|\tuET\|_{L^4(\Omega)}^{8/3}}{2\nu \E(\tuET)}
\label{eq:Theta}
\end{equation}
where we used the energy equation \eqref{eq:dKdt}.  This quantity is
also plotted in Figure \ref{fig:P3vsK}b where we see that for each
value of $\K_0$ we have $\Theta_{\K_0} > \Xi_{\K_0}(T)$, $T > 0$. Its
dependence on $\K_0$ is approximated by the power-law relation
\begin{equation}
\Theta_{\K_0}  \approx (4.963\pm 0.492) \K_0^{0.32\pm 0.03}
\label{eq:ThetaK0}
\end{equation}
from which we deduce that the quantities $\Theta_{\K_0}$ and
$\Xi_{\K_0}(T)$ exhibit a weaker growth with $\K_0$ than given by the
expression on the RHS in \eqref{eq:Xi} where the exponent is 2/3.
This thus demonstrates that estimate \eqref{eq:LPSbound} is not sharp
for any time window $T$.}

\FloatBarrier

\section{Discussion and Conclusions}
\label{sec:final}

In this study we have undertaken a systematic computational search for
potential finite-time singularities in incompressible Navier-Stokes
flows based on the Ladyzhenskaya-Prodi-Serrin conditional regularity
criterion \eqref{eq:LPS}. This criterion asserts that a solution
$\u(t)$ is smooth and satisfies the Navier-Stokes system \eqref{eq:NS}
in the classical sense on the time interval $[0,T]$ provided the
integral $\int_0^T \| \u(\tau) \|_{L^q(\Omega)}^{4q/(3(q-2))} \,
d\tau$, where $q > 3$, is bounded. In our study we chose $q = 4$ and
$p = 8$ which is the pair of integer-valued indices closest to the
critical case with $p = 3$. To the best of our knowledge, this is the
first such investigations based on the Ladyzhenskaya-Prodi-Serrin
condition \eqref{eq:LPS} and it complements earlier studies based on
the enstrophy condition \cite{ld08,ap16,KangYumProtas2020}.

The idea of our approach is to consider classical solutions of the
Navier-Stokes system \eqref{eq:NS} which might blow up in finite time.
Initial data which might potentially lead to a singularity is sought
by solving Problems \ref{pb:PhiL4} and \ref{pb:PhiH3/4} in which
quantity \eqref{eq:Phi} is maximized subject to different sets of
constraints. These problems were solved numerically with a
state-of-the-art adjoint-based maximization approach formulated in the
continuous (infinite-dimensional) setting.  Since such approaches are
most conveniently defined in Hilbert spaces, our optimal initial data
was sought in the space $H^{3/4}(\Omega)$, which is the largest
Sobolev space with Hilbert structure embedded in the space
$L^{4}(\Omega)$ appearing in condition \eqref{eq:LPS} when $q = 4$.

Problems \ref{pb:PhiL4} and \ref{pb:PhiH3/4} both admit two branches
of maximizing solutions for a broad range of constraint values,
cf.~Figures \ref{fig:maxLPSvsT}a and \ref{fig:maxLPSvsT}b. {It is
  interesting to note that while Problems \ref{pb:PhiL4} and
  \ref{pb:PhiH3/4} involve the same objective functional
  $\Phi_T(\u_0)$ maximized over the same function space
  $H^{3/4}(\Omega)$, but subject to different, though related,
  constraints, their solutions are in fact very different.  However,
  in none of the cases was there any evidence found} for emergence of
a singularity, in the sense that quantity \eqref{eq:Phi} remains
bounded for all values of the constraints and all optimization windows
$T$. However, when considering the corresponding growth of enstrophy,
solutions of Problem \ref{pb:PhiL4} from the partially-symmetric
branch were found to attain enstrophy values scaling in proportion to
$\E_0^{3/2}$, cf.~\eqref{eq:fit_maxET_L}.  This is interesting because
the same power-law dependence (but with a different, larger,
prefactor) of the maximum attained enstrophy on $\E_0$ was obtained in
Navier-Stokes flows with initial data constructed to maximize the
finite-time growth of enstrophy in \cite{KangYumProtas2020},
cf.~Figure \ref{fig:maxE}, as well as in 1D Burgers flows with initial
data determined in an analogous manner \cite{ap11a}. Thus, extreme
Navier-Stokes flows with distinct structure obtained by maximizing two
different quantities are characterized by the same power-law relation
$\max_{t \ge 0} \E(t) \sim \E_0^{3/2}$ describing the dependence of
the maximum attained enstrophy on the initial enstrophy. We recall
that at present there are no rigorous a priori bounds on the growth of
enstrophy and the best available estimate \eqref{eq:Et_estimate_E0}
has an upper bound which becomes infinite in finite time.

As the second main contribution of our study, we have considered the a
priori estimate \eqref{eq:LPSbound} and showed that it does not appear
sharp, although the degree to which the expression on the RHS
overestimates the growth of $\frac{1}{T} \int_0^T \| \u(\tau)
\|_{L^4(\Omega)}^{8/3} \, d\tau$ with $\K_0$ is reduced as $T
\rightarrow \infty$ (by ``sharpness'' we mean that the expression on
the LHS in the estimate scales with $\K_0$ in the same way up to a
prefactor as the upper bound on the RHS).  This observation was
deduced by solving Problem \ref{pb:PsiL2} for a range of values of
$\K_0$ and $T$, and then extrapolating the results to large values of
$T$. This lack of sharpness appears to be a consequence of the fact
that the term $\|\u(T)\|_{L^2(\Omega)}^2$, which is dropped in
\eqref{eq:LPSbound2a}, is in general non-negligible for finite $T$,
but becomes less significant as $T \rightarrow \infty$. These results
thus demonstrate that estimate \eqref{eq:LPSbound} may potentially be
improved by reducing the power of $\K_0$ in the upper bound. This
should not come as a surprise since the instantaneous estimate
\eqref{eq:dLqdt} was found not to be sharp as well, cf.~Figure
\ref{fig:dLqdt} and relation \eqref{eq:fit_dLqdt}.  We emphasize,
however, that given the fact that Problems \ref{pb:PhiL4},
\ref{pb:PhiH3/4} and \ref{pb:PsiL2} are non-convex, the observations
made above cannot be regarded as definitive, since it is possible that
despite our efforts we might not have found global maximizers.

As regards future studies, it is worthwhile to reconsider the problems
investigated here using a formulation where the optimal initial data
is sought directly in the space $L^4(\Omega)$ rather than in
$H^{3/4}(\Omega)$. This can be done using an extension of the
adjoint-based optimization approach we used to more general Banach
spaces \cite{protas2008}, which is however more technically involved.
It is also interesting to probe the Ladyzhenskaya-Prodi-Serrin
criterion \eqref{eq:LPS} for a broad range of values of $p$ and $q$,
as well as to consider generalizations of this criterion involving
derivatives of different order of the velocity field obtained in
\cite{Gibbon2018}. The limiting (critical) case with $q = 3$,
cf.~\eqref{eq:LPS3}, is particularly interesting. However, given the
non-differentiability of the norm $\|\cdot\|_{L^\infty([0,T])}$, this
problem is not amenable to straightforward solution with the
gradient-based optimization approach considered here. On the other
hand, condition \eqref{eq:LPS3} can be probed by maximizing the
finite-time growth of the norm $\|\u(T)\|_{L^3(\Omega)}$, in analogy
to Problem \ref{pb:maxET} solved in \cite{KangYumProtas2020}. Finally,
it is also of interest to consider the problems studied here on the
unbounded domain $\RR^3$ instead of a torus.

\section*{Acknowledgments}

This work is dedicated to the memory of the late Charlie Doering, our
dear friend and collaborator, who inspired us to pursue this research
direction.  The authors wish to express thanks to John Gibbon, Evan
Miller and Koji Ohkitani for enlightening and enjoyable discussions.
They also acknowledge the support through an NSERC (Canada) Discovery
Grant.  Computational resources were provided by Compute Canada under
its Resource Allocation Competition.

\appendix


\section{Derivation of Estimate \eqref{eq:LPSbound} with an
  Explicit Upper Bound}
\label{sec:LPSbound}

We begin with the Gagliardo--Nirenberg inequality
\begin{equation}
||D^{j}\u||_{L^{p}}\leq C||D^{m}\u||_{L^{r}}^{\alpha} \, ||\u||_{L^{q}}^{1-\alpha}, \quad 
\text{where} \
\frac{1}{p}=\frac{j}{n}+\left(\frac{1}{r}-\frac{m}{n}\right)\alpha+\frac{1-\alpha}{q}
\quad \text{and} \quad 
\frac{j}{m}\leq \alpha \leq 1.
\label{eq:GN}
\end{equation}
Setting $j=0$, $m=1$, $r=2$, $q=2$, and $n=3$ we obtain
$\frac{1}{p}=0+\left(\frac{1}{2}-\frac{1}{3}\right)\alpha+\frac{1-\alpha}{2}$
and $\alpha=\frac{3(p-2)}{2p}$, such that for $2\leq p\leq6$
inequality \eqref{eq:GN} becomes
\begin{equation}
||\u||_{L^{p}} \leq C||\bnabla \u||_{L^{2}}^{\alpha} \, ||\u||_{L^{2}}^{1-\alpha}.
\label{eq:GN0}
\end{equation}
Raising both sides of \eqref{eq:GN0} to the power $2/\alpha$,
integrating with respect to time over $[0,T]$ and then using the
energy equation \eqref{eq:dKdt} yields
\begin{subequations}
\label{eq:LPSbound2}
\begin{align}
\int_{0}^{T}||\u(t)||_{L^p}^{\frac{4p}{3(p-2)}} \, dt & \leq  \int_{0}^{T}C ||\bnabla \u(t)||_{L^{2}}^{2} \, ||\u(t)||_{L^{2}}^{\frac{2(1-\alpha)}{\alpha}} \, dt \nonumber \\
 & \leq  C\left(||\u_{0}||_{L^{2}}\right)^{\frac{2(1-\alpha)}{\alpha}}\int_{0}^{T}||\bnabla \u(t) ||_{L^{2}}^{2} \, dt \nonumber \\
 & =  \frac{C}{2\nu}||\u_{0}||_{L^{2}}^{\frac{2(1-\alpha)}{\alpha}}\left(||\u_{0}||_{L^{2}}^{2}-||\u(T)||_{L^{2}}^{2}\right)  \label{eq:LPSbound2a} \\
 & \leq  C \K_0^{\frac{2p}{3(p-2)}}, \qquad 2\leq p\leq 6.
\label{eq:LPSbound2b}
\end{align}
\end{subequations}
On the other hand, we can deduce from \cite[Theorem 2(i)]{Gibbon2018}
that
\begin{equation}
\int_{0}^{T}||\u(t)||_{L^{p}}^{\frac{p}{p-3}} \, dt \leq  
C\left(\int_{0}^{T}||\u(t)||_{L^{2}}^{2} \, dt\right)^{\frac{3}{2}} \leq 
C\K_0^{3}, \qquad p>6.
\label{eq:LPSbound6}
\end{equation}
The ranges of validity of estimates \eqref{eq:LPSbound2} and
\eqref{eq:LPSbound6} do not overlap, however, in the borderline case
when $p=6$ the expressions on the LHS in the two estimates coincide,
yet the upper bound in the first estimate is $C \K_0$ and therefore
has a smaller exponent than the upper bound in the second estimate.

\section{Problem \ref{pb:PhiL4} in the Limit $T \rightarrow 0$}
\label{sec:dLqdt}

In this Appendix we show that solutions of Problem \ref{pb:PhiL4}
approximate solutions of the instantaneous optimization problem
$\max_{\u \in  {\L}_{B}} \frac{d}{dt} \|\u\|_{L^q(\Omega)}^q$, $q > 3$ in
the limit $T \rightarrow 0$. We have for $p,q \ge  1$
\begin{align*}
 \frac{d \|\u(t)\|_{L^q(\Omega)}^p}{dt}\bigg|_{t=0} & = \frac{\|\u(T)\|_{L^q(\Omega)}^p - \|\u_0\|_{L^q(\Omega)}^p}{T} + \O(T) \\
& = \frac{\frac{d}{dT}\int_0^T \|\u(t)\|_{L^q(\Omega)}^p \, dt - \|\u_0\|_{L^q(\Omega)}^p}{T} + \O(T) \\
& = \frac{\frac{1}{T}\int_0^T \|\u(t)\|_{L^q(\Omega)}^p \, dt + \O(T) - \|\u_0\|_{L^q(\Omega)}^p}{T} + \O(T),
\end{align*}
where we used the first-order finite-difference approximation of the
derivative twice and the fundamental theorem of calculus.  Then, after
taking the maximum on both sides we obtain for $T \rightarrow 0$
\begin{equation}
\max_{\u \in  {\L}_{B}} \frac{d \|\u(t)\|_{L^q(\Omega)}^p}{dt}\bigg|_{t=0} = 
\frac{\frac{1}{T}\max_{\u \in  {\L}_{B}}  \int_0^T \|\u(t)\|_{L^q(\Omega)}^p \, dt - \|\u_0\|_{L^q(\Omega)}^p}{T} + \O(1).
\label{eq:maxdLqdt}
\end{equation}
Finally, to be able to relate this result to estimate
\eqref{eq:dLqdt}, we apply the chain rule to obtain
\begin{equation}
\max_{\u \in  {\L}_{B}} \frac{d \|\u(t)\|_{L^q(\Omega)}^q}{dt}\bigg|_{t=0} = \frac{q}{p} \|\u_0\|_{L^q(\Omega)}^{q-p}\left(
\frac{\frac{1}{T}\max_{\u \in  {\L}_{B}}  \int_0^T \|\u(t)\|_{L^q(\Omega)}^p \, dt - \|\u_0\|_{L^q(\Omega)}^p}{T} + \O(1)\right).
\end{equation}


\end{document}